\documentclass[lettersize,journal]{IEEEtran}
\usepackage{amsmath,amsfonts}
\usepackage{algorithmic}
\usepackage{algorithm}
\usepackage{array}
\usepackage[caption=false,font=normalsize,labelfont=sf,textfont=sf]{subfig}
\usepackage{textcomp}
\usepackage{stfloats}
\usepackage{url}
\usepackage{verbatim}
\usepackage{graphicx}
\usepackage{cite}

\usepackage{amsmath}
\usepackage{tikz}
\usepackage{mathdots}
\usepackage{yhmath}
\usepackage{cancel}
\usepackage{color}
\usepackage{array}
\usepackage{multirow}
\usepackage{amssymb}
\usepackage{tabularx}
\usepackage{extarrows}
\usepackage{booktabs}
\usepackage{cellspace}
\usepackage{hyperref}
\usepackage{cleveref}
\usepackage{float}
\usepackage{makecell}
\usepackage{tabularx}
\usepackage{longtable}
\usepackage{bm}

\usetikzlibrary{patterns}
\usetikzlibrary{shadows.blur}
\usetikzlibrary{shapes}

\begin{document}

\title{Learning Noise-Robust Stable Koopman Operator for Control with Hankel DMD}

\author{Shahriar Akbar Sakib, Shaowu Pan
\thanks{Manuscript submitted on August 13, 2024. This work was supported in part by the startup grant from Rensselaer Polytechnic Institute. The code will be made available upon publication at \texttt{https://github.com/csml-rpi/stirk}.}
\thanks{S. A. Sakib is with the Department of Mechanical, Aerospace, and Nuclear Engineering at Rensselaer Polytechnic Institute, Troy, NY 12180 USA. (e-mail:sakibs@rpi.edu)}
\thanks{S. Pan is in the Department of Mechanical, Aerospace and Nuclear Engineering at Rensselaer Polytechnic Institute, Troy, NY 12180 USA. (e-mail:pans2@rpi.edu)}
}

\markboth{Journal of \LaTeX\ Class Files,~Vol.~14, No.~8, August~2021}%
{Shell \MakeLowercase{\textit{et al.}}: A Sample Article Using IEEEtran.cls for IEEE Journals}

\IEEEpubidadjcol


\maketitle

\begin{abstract}
We propose a noise-robust learning framework for the Koopman operator of nonlinear dynamical systems, with guaranteed long-term stability and improved model performance for better model-based predictive control tasks. Unlike some existing approaches that rely on ad hoc observables or black-box neural networks in extended dynamic mode decomposition (EDMD), our framework leverages observables generated by the system dynamics, when the system dynamics is known, through a Hankel matrix, which shares similarities with discrete Polyflow. When system dynamics is unknown, we approximate them with a neural network while maintaining structural similarities to discrete Polyflow. To enhance noise robustness and ensure long-term stability, we developed a stable parameterization of the Koopman operator, along with a progressive learning strategy for rollout loss. To further improve the performance of the model in the phase space, a simple iterative data augmentation strategy was developed. Numerical experiments of prediction and control of classic nonlinear systems with ablation study showed the effectiveness of the proposed techniques over several state-of-the-art practices. 
\end{abstract}


\section{Introduction}
\label{sec:introduction}
Data-driven methods have shown promising potential in modeling and control of non-linear systems~\cite{brunton2021modern,van2020data}.
The Koopman operator, as originally proposed by Bernard Koopman~\cite{koopman1931hamiltonian}, plays a pivotal role by providing a global linear representation of nonlinear dynamical systems~\cite{mezic2005spectral,budivsic2012applied}.
Interested readers may refer to several excellent reviews on the Koopman operator~\cite{brunton2021modern, mezic2013analysis,otto2021koopman,bevanda2021koopman,mezic2023operator}. 
The Koopman operator has been applied in various fields, including robust prediction of COVID-19~\cite{mezic2024koopman}, reinforcement learning~\cite{rozwood2024koopman}, soft robotics control~\cite{bruder2021koopman}, stability analysis of genetic toggle switches~\cite{harrison2021stability}, and modeling turbulent shear flows~\cite{constante2024data}.

The biggest advantage of the Koopman operator for the dynamics and control community is its ability to facilitate the use of well-established linear control/observer techniques for nonlinear systems~\cite{korda2016linear,kaiser2021data,surana2016linear}, provided that a good approximation of the global linear representation is learned from the data. 
For example, model predictive control (MPC) is typically computationally expensive as it requires solving a non-convex optimization problem online.
Traditional linearization techniques that approximate the nonlinear system around a local reference point might result in suboptimal MPC performance and could struggle with significant deviations from this point. 
On the other hand, the linear representation of the Koopman operator enables the transformation of the nonconvex optimization problem in nonlinear MPC into a convex optimization problem in linear MPC, as demonstrated by \cite{korda2018linear, peitz2020data}.
Beyond MPC, the linear representation provided by the Koopman operator can also benefit LQR~\cite{kaiser2021data}, anomaly detection~\cite{qian2021koopman}, and fault-tolerant control~\cite{bakhtiaridoust2024koopman} for nonlinear systems.  

Despite its successful applications, several challenges remain in utilizing the Koopman operator.
The first challenge is to find an appropriate invariant subspace for the Koopman operator. This involves selecting suitable basis functions that span an invariant subspace. Common choices of basis functions include radial basis functions (RBFs)~\cite{li2017extended} or polynomial functions~\cite{brunton2016koopman}.  
Interestingly, Otto et al. \cite{otto2024learning} proposed learning the latent linear model and hidden state jointly, thereby eliminating the need to explicitly specify a dictionary of functions, as required by existing methods such as extended dynamic mode decomposition (EDMD). However, this approach requires solving a nontrivial optimization problem for the linear dynamics and hidden states, which are unique to each trajectory.

Recently, the concept of polyflow was introduced by \cite{jungers2019non}. Polyflow employs the governing equations of the system to construct basis functions, offering a non-local linearization approach for nonlinear differential equations. This concept is closely related to the derivative-based approximation of Koopman operator~\cite{mamakoukas2021derivative} except closure is considered.  In earlier work, it was observed by Mezi\'{c} \cite{mezic2022numerical} that instead of relying on arbitrary sets of observables, one could use observables generated by the dynamics, such as time delays that form a Krylov subspace, to study the spectral properties of the Koopman operator. This approach has also been applied in the context of Hankel-DMD methods.
Despite the initial success of Polyflow for MPC~\cite{wang2021immersion}, merely minimizing the one-step-ahead error remains susceptible to noise in the data.
The second challenge is to ensure long-term prediction accuracy. 
Wu et al.~\cite{wu2021challenges} revealed that even under mild noise and nonlinearity, DMD struggles to accurately predict system dynamics. Several recent approaches have been proposed to address the challenges posed by measurement noise in learning Koopman operators. Sinha et al. \cite{sinha2023online} developed a robust online Koopman framework that employs a min-max optimization strategy to mitigate the effect of noisy streaming data, while Hao et al. \cite{hao2024deep} proposed a deep Koopman learning method that incorporates considerations of bounded noise by updating observable functions during training. However, these methods either require restrictive assumptions on noise bounds or face challenges in ensuring long-term prediction stability, which motivates the need for our stable parameterization and rollout-based learning strategy. 
Although the use of roll-out techniques for recurrent loss has proven beneficial when data are limited and noisy both in the deep learning for the Koopman operator~\cite{otto2019linearly,lusch2018deep,pan2020physics} and the reduced-order modeling community \cite{uy2023operator}, optimizing recurrent loss can be tricky and prone to poor local minimum in practice.  
The third challenge in Koopman operator-based control is to ensure the stability of the Koopman operator itself, since the real part of any continuous-time eigenvalue of the Koopman operator of any physical system should not be positive~\cite{pan2020physics,bevanda2022diffeomorphically,bevanda2022learning}. If the forward pass becomes unstable, it inevitably destabilizes the backward pass as well. 
Lastly, curating training data for Koopman operator-based control often depends on ad-hoc selections of input signals and initial conditions, which can result in insufficient sampling especially in the high-dimensional phase space and, consequently, poor control performance.
Besides, it is important to note that we do not attempt to address other challenges such as the convergence of Koopman operator. Interested readers are referred to this excellent paper by Colbrook and Mezi\'{c}~\cite{colbrook2024limits}. 
%

In order to address the above issues, the goal of this paper is to develop a noise-robust learning framework for the Koopman operator. Specifically, the contributions of this work are as follows:
\begin{itemize}
    \item improving noise-robustness and long-term prediction accuracy via minimizing recurrent loss via roll-out,
    \item stabilizing Koopman Operator by negative-definite parameterization,
    \item iteratively learning Koopman-based model on augmented data with previous responses of the controlled nonlinear system.
\end{itemize}

The structure of the paper is as follows: \Cref{sec:PROPOSED FRAMEWORK} describes the proposed framework, including the use of the polyflow basis and the stabilized formulation.
\Cref{sec:Numerical Examples} presents the results of our experiments, demonstrating the efficacy of our approach. 
\Cref{sec:adaptive data augmentation} discusses the results of applying iterative data augmentation on Koopman-based control tasks. 
In \Cref{sec:ablation}, we performed an ablation study to evaluate the contributions of the various methods used in this paper. 
\Cref{sec:Conclusions} provides the conclusion of the paper.


\section{Proposed Framework}
\label{sec:PROPOSED FRAMEWORK}
\subsection{Koopman Operator for Controlled System}
To provide a rigorous justification for constructing a linear predictor by lifting the state space, we briefly outline the Koopman operator approach to analyze a controlled dynamical system, following the approach of \cite{korda2016linear}. 
We consider the following discrete-time nonlinear system with control,
\begin{align}
x^+ & =f(x,u),
\label{eq:nonlinear_system}
\end{align}
where $x \in \mathbb{R}^{n}$ is the state of the system, $x^{+}$ is the system state at the next time step, and $u \in \mathcal{U} \subset \mathbb{R}^{m}$ is the control input. $\mathcal{U}$ is the space of the available control inputs. $f:\mathbb{R}^{n} \times \mathbb{R}^{m}\rightarrow \mathbb{R}^{n}$ is the transition mapping. $x^+$ is the state of the system for the next time step. 
We define the extended state comprising the state $x$ of the system and the control sequence $\mathbf{u} \in \ell(\mathcal{U})$, where $\ell(\mathcal{U})$ is the space of all sequences $(u_i)_{i=0}^{\infty}$ with $(u_i)_{i=0}^{\infty} \in \mathcal{U}$,
\[
\chi = \begin{bmatrix}
x \\ \mathbf{u}
\end{bmatrix},
\]
The dynamics of the extended state $\chi$ is described by
\begin{align}
\chi^+ = F(\chi) := \begin{bmatrix}
f(x, \mathbf{u}(0)) \\
\mathcal{S}\mathbf{u}
\end{bmatrix},    
\label{eq:extended_dynamics}
\end{align}
where $\mathcal{S}$ is the left shift operator, i.e., $(\mathcal{S}\mathbf{u})(i) = \mathbf{u}(i + 1)$, and $u(i)$ denotes the $i$th element of the sequence $\mathbf{u}$.
The Koopman operator $\mathcal{K}: \mathcal{H} \to \mathcal{H}$ associated with \cref{eq:extended_dynamics} is defined by
\[
(\mathcal{K}\psi)(\chi) = \psi(F(\chi))
\]
for any observable $\psi : \mathbb{R}^n \times \ell(\mathcal{U}) \to \mathbb{R}$ belonging to some space of observables $\mathcal{H}$.
Despite the fact that the Koopman operator $\mathcal{K}$ is a linear operator that fully describes the time evolution of any observable function in a non-linear dynamical system, it is inherently infinite-dimensional which makes learning the Koopman operator challenging. 
In the following subsection, we adopt the widely used concept of \emph{lifted linear predictor}~\cite{bevanda2023koopman,bevanda2024kernel} as a finite-dimensional approximation of the Koopman operator. 

%

\subsection{Lifted Linear Predictor}
Let us use a subscript to indicate the time step and rewrite \cref{eq:nonlinear_system} as follows. 
\begin{equation}
    x_{k+1} = f(x_k, u_k),
    \label{eq:discrete-time-nonlinear-dynamical-system}
\end{equation}
where \( f: \mathbb{R}^n \times \mathbb{R}^m \to \mathbb{R}^n \) represents a general non-linear mapping, and \( x_k = x(t_k) \in \mathbb{R}^n \) and \( u_k = u(t_k) \in \mathbb{R}^m \) denote the state of the system and the input of the control at time \( t_k \), respectively. Given an initial condition \( x_0 \) at time \( t_0 \) and a sequence of control inputs \( \{u_k\}_{k=0}^{N_p} \), where \( N_p \) is the length of the prediction horizon, the objective is to estimate the sequence of future states \( \{ x_k \}_{k=1}^{N_p+1} \). 

To approximate the evolution of this nonlinear system using a linear framework, following the approach of \cite{korda2016linear}, we consider a \textit{lifted linear predictor}, modeled as a linear dynamical system parameterized by a triplet of matrices \( (A, B, C) \):
\begin{equation}
\begin{aligned}
    z_0 &= \Phi(x_0), &  \\
    z_{k+1} &= A z_k + B u_k, & k = 0, 1, \dots, N_p, \\
    \hat{x}_k &= C z_k, & k = 1, \dots, N_p+1.
\end{aligned}
\label{eq:lifted-linear-predictor}
\end{equation}
where \( z_k \in \mathbb{R}^{N_\phi} \) is the \textit{lifted state} in an extended space of dimension \( N_\phi \), and \( \Phi: \mathbb{R}^n \to \mathbb{R}^{N_\phi} \) is a non-linear \textit{ lifting function} that maps the original state space to a higher dimensional feature space. The vector \( \hat{x}_k \in \mathbb{R}^n \) represents the estimated state of the system at time $t_k$.

The lifting function \( \Phi(\cdot) \) is defined as:
\begin{equation}
    \Phi(x) = \begin{bmatrix} \phi_1(x) \\ \phi_2(x) \\ \vdots \\ \phi_{N_\phi}(x) \end{bmatrix} \in \mathbb{R}^{N_\phi},
\end{equation}
where each function \( \phi_i: \mathbb{R}^n \to \mathbb{R} \) is an \textit{observable}, capturing potentially non-linear transformations of the state. The choice of observables \( \{\phi_i(x)\}_{i=1}^{N_\phi} \) determines the accuracy of the linear approximation.

It is important to note that, under certain conditions (e.g., control-affine systems), one can also adopt a \emph{bilinear} model realization \cite{zhao2024deep, bruder2021advantages, yu2022autonomous, goswami2021bilinearization, folkestad2022koopnet}. Suppose that there exist coefficient matrices
$
A \in \mathbb{R}^{N_\phi \times N_\phi}, 
\quad
B_i \in \mathbb{R}^{N_\phi \times N_\phi}, 
\quad
C \in \mathbb{R}^{n \times N_\phi},$
such that for
$
z_k = [z_{k}^{(1)}, \dots, z_{k}^{(N_\phi)}]^\top, 
\quad
x_k = [x_{k}^{(1)}, \dots, x_{k}^{(n)}]^\top, 
\quad
u_k = [u_{k}^{(1)}, \dots, u_{k}^{(m)}]^\top,
$
the dynamics takes the form
\begin{equation}
\begin{aligned}
z_{k+1} &= Az_k +\sum _{i=1}^{m} u_{k}^{( i)} B_{i} z_{k},\\
x_k &= C \, z_k,
\end{aligned}
\label{eq:bilinear_definition}
\end{equation}
then $\bigl(A,\{B_i\}_{i=1}^{m},C\bigr)$ is said to be a bilinear realization of the system. Although extending our approach to the bilinear setting is natural, here we opt for a linear formulation in \cref{eq:lifted-linear-predictor} for simplicity and seamless integration with standard control methods~\cite{shang2024willems}.

\subsection{Extended DMD for Control}
For practical purposes, we adapt the Extended Dynamic Mode Decomposition (EDMD)~\cite{williams2015data} approach followed by Korda and Mezi\'{c}~\cite{korda2018linear} for control. 
The EDMD provides a straightforward method to approximate the Koopman operator from the data, thereby obtaining the linear predictor lifted in \cref{eq:lifted-linear-predictor}. Since our goal is not to predict the control sequence, the observable function will depend only on the state of the system. Thus, the process begins by selecting \( N_\phi \) observables \( \phi_1, \ldots, \phi_{N_\phi} \) and forming the lifted state \( z \) from the measured state \( x \) of the system \eqref{eq:discrete-time-nonlinear-dynamical-system} as follows:
\begin{equation}
z = \Phi(x) = [ \phi_1(x), \phi_2(x), \ldots, \phi_{N_\phi}(x) ]^\top.
\end{equation}

The one-step prediction loss based EDMD algorithm assumes that a collection of data \((x_j, u_j, x_j^+), j = 1,2, \ldots, K\) satisfying \(x_j^+ = f(x_j,u_j)\) is available and seeks matrices $A$ and $B$ minimizing
\begin{equation}
    \sum_{j=1}^{K} \|\Phi(x_j^+) - A\Phi(x_j) - Bu_j\|_2^2,
    \label{eq:one-step loss}
\end{equation}

\subsection{Polyflow Basis Functions for EDMD}
To construct suitable basis functions \( \Phi \) that span an invariant subspace for the Koopman operator, we leverage \textit{Polyflow basis functions}, inspired by the work of Wang and Jungers~\cite{wang2021immersion}, as well as Krylov subspace methods discussed by Mezi\'{c}~\cite{mezic2022numerical}. These basis functions provide a physics-informed approach to lifting the system state, ensuring a structured representation of the non-linear dynamics.

\subsubsection{Known Dynamics}
When the dynamics of the system \( f(x, u) \) is explicitly known, we define the Polyflow basis functions of order \( N \), denoted as \( T_N(x) \), as follows:
\begin{equation}
\label{eq:polyflow}
T_{N}(x) :=
\begin{bmatrix}
    x \\
    f(x,0) \\
    f(f(x,0),0) \\
    \vdots \\
    f^{(N-1)}(x,0)
\end{bmatrix} \in \mathbb{R}^{N_\phi}.
\end{equation}
Here, \( f^{(i)}(x,0) \) represents the \( i \)-th recursive application of \( f \) with zero input, which captures the intrinsic evolution of the system. The dimension of the lifted state is $N_{\phi} = N\times n$.

Using these Polyflow basis functions to lift the state of the system, we obtain the following linear system representation, consistent with the linear predictor lifted in \cref{eq:lifted-linear-predictor}:
\begin{equation}
\begin{aligned}
    z_0 &= T_N(x_0), \\
    z_{k+1} &= A z_k + B u_k, \\
    \hat{x}_k &= C z_k.
\end{aligned}
\label{eq:linear_embedding}
\end{equation}
The matrix \( C \) is typically determined by the selection of \( T_N(x) \) and encodes the mapping from the lifted state space back to the original state space.

\subsubsection{Unknown Dynamics}
If the governing dynamics \( f(x, u) \) is unknown, we approximate them using a neural network \( g{_\theta}(x) : \mathbb{R}^n \mapsto \mathbb{R}^n \) parameterized by \( \theta \), while maintaining the recursive structure of Polyflow. The corresponding basis functions are then defined as:
\begin{equation}
\label{eq:unknown_recursive}
T_{N}(x) :=
\begin{bmatrix}
    x \\
    g{_\theta}(x) \\
    g_{\theta}(g{_\theta}(x)) \\
    \vdots \\
    g_{\theta}^{(N-1)}(x)
\end{bmatrix} \in \mathbb{R}^{N_\phi}.
\end{equation}
Here, \( g_{\theta}(x) \) serves as a learned approximation of \( f(x,0) \), recursively applied to construct the lifted state representation. The dimension of the lifted state is $N_{\phi} = N\times n$. This data-driven extension allows Polyflow to be applied even in scenarios where an explicit model of the system dynamics is unavailable.

\subsubsection{Connection to Hankel DMD}
Our proposed use of Polyflow basis functions bears conceptual similarity to the Hankel DMD method~\cite{arbabi2017ergodic,brunton2017chaos}, as both leverage delay-embedded observables to approximate the action of the Koopman operator through a Krylov subspace structure. Specifically, the Polyflow recursion we employ mimics delay embeddings by iteratively composing either the true or learned dynamics to generate basis functions. However, unlike Hankel DMD~\cite{brunton2017chaos}---which constructs \emph{purely numerical} delay-coordinate snapshots from raw trajectory data---Polyflow builds its \emph{function} basis through recursive application of (possibly learned) dynamics, resulting in a physics-informed lifting map in an EDMD fashion. 
This makes Polyflow particularly advantageous when historical data is unavailable for starting.
%
%

\subsection{Roll-out Loss}

To improve long-term prediction accuracy and noise resistance in training data, we employ a \textit{multi-step roll-out loss} instead of the one-step prediction loss \cref{eq:one-step loss}. The multi-step roll-out loss refers to cumulative prediction error in multiple future time steps (also known as rollout length). The rollout length $R$ determines how many steps in the future the model is trained to predict.

Consider a data set consisting of \( M \) trajectories, each trajectories consisting of \( R \)  roll-out steps,
\( \{(x_{0,j}, x_{1,j}, \dots, x_{R,j}, u_{0,j}, \dots, u_{R-1,j})\}_{j=1}^{M} \), 
where \( x_{i,j} \) denotes the state of the system at the \( i \) -th time step of the \( j \)-th trajectory, and \( u_{i,j} \) is the corresponding control input.

The optimization problem for learning the system matrices \( A \) and \( B \) is then formulated as follows:
\begin{equation}
\underset{A,B}{\min} \sum_{j=1}^{M} \sum_{r=1}^{R} 
\left\| T_N(x_{r,j}) - A^r T_N(x_{0,j}) - \sum_{k=0}^{r-1} A^{r-k-1} B u_{k,j} \right\|_2^2,
\label{eq:optimization_with_roll-out}
\end{equation}
where \( A \) and \( B \) are the unknown system matrices to be learned, and \( T_N(x) \) represents the transformation of the Polyflow basis function of the state of the system.

The hyperparameters \( N \) and \( R \) are typically determined by tuning the hyperparameters to balance the complexity of the model and predictive performance. Adjusting these parameters is a trade-off: increasing
$N$ can capture more complex dynamics but can lead to overfitting if not managed carefully, while higher values of $R$ typically enhance long-term prediction accuracy at the cost of increased computational burden. 

\subsection{Curriculum Learning}
Curriculum learning~\cite{bengio2009curriculum} is a training strategy in which a model is first trained in simpler tasks and then gradually exposed to more challenging ones. This approach helps the model to learn and adapt progressively to increasing complexity. 

In our framework, we implement a form of curriculum learning by progressively increasing the length of the rollout \( R \) during the training process. \emph{Progressive scheduling} refers to gradually increasing the complexity of the training task - in our case, by increasing the length of roll-out
$R$ during training. This aligns with the principles of curriculum learning, where the model first learns easier subtasks before tackling harder ones. We refer to this approach as \emph{progressive rollout schedule}.

Initially, the model is trained with shorter rollout lengths, focusing on immediate predictions. As training progresses and the model stabilizes, the roll-out length \( R \) increases geometrically after certain epochs. This gradual increase is expected to allow the model to learn from easier short-term predictions before tackling the more complex task of long-term prediction, thus aligning with the curriculum learning paradigm. In doing so, the model is expected to build its capacity to handle extended sequences, improving overall prediction accuracy over long time horizons. Compared to direct exposure to long rollout lengths from the start, this progressive approach offers optimization benefits by improving stability and convergence. Training with long roll-outs early on can lead to high variance in gradients and poor local optima due to compounding prediction errors. By starting with shorter roll-outs, the model establishes a more robust foundation before incrementally handling extended sequences, ultimately enhancing prediction accuracy over long time horizons.

\subsection{Stabilized Formulation}
\label{sec:stabilized_formulation}

To guarantee the stability of both the Koopman matrix $A$ and the training process under rollout loss, we enforce the eigenvalues of $A$ to be stable, employing the fact that a skew-symmetric matrix with negative real-valued diagonal elements has stable or neutral stable eigenvalues in continuous-time dynamical system~\cite{pan2020physics}. This guarantees the stability of the Koopman matrix during training with roll-out loss, preventing instability in both forward and backward passes.

For some matrix \( Q \in \mathbb{R}^{\tilde{n} \times \tilde{n}} \), $\sigma_1, ..., \sigma_{\tilde{n}} \in \mathbb{R}$, a diagonal matrix \( C =  \mathrm{diag}(-\sigma_1^2,...,-\sigma_{\tilde{n}}^2)\in \mathbb{R}^{\tilde{n} \times \tilde{n}} \) with negative real diagonal elements, we construct:
\begin{equation}
    \tilde{A} = Q - Q^{\top} + C.
\end{equation}
By construction, the eigenvalues of \( \tilde{A} \) lie in the left half of the complex plane, ensuring stability for the continuous-time dynamical system.

To obtain the discrete-time matrix \( \hat{A} \in \mathbb{R}^{\tilde{n} \times \tilde{n}} \), we apply the matrix exponential to \( \tilde{A} \) scaled by the sampling timestep \( \Delta t \):
\begin{equation}
    \hat{A} = \textrm{expm}\left( \tilde{A} \Delta t \right),
\end{equation}
where expm refers to matrix exponential. 
Since the Koopman operator is typically expressed in a different coordinate system, a similarity transformation using a coordinate transformation matrix \( P \in \mathbb{R}^{\tilde{n} \times \tilde{n}} \) is required to recover \( A \):
\begin{equation}
    A = P \hat{A} P^{-1}.
\end{equation}

This formulation ensures that the Koopman operator remains stable throughout the training, which we refer to as the \textit{ dissipative Koopman operator}. In contrast, the \textit{ standard Koopman operator} refers to a lifted system where \( A \) is unconstrained.

\subsection{Model Predictive Control}
One of the most important applications of the lifted linear predictor is to provide a surrogate model for MPC, which is known as Koopman-MPC (KMPC)~\cite{korda2018linear}. At each discrete time step, the MPC solves the following optimization problem on the prediction horizon \( N_p \):
\begin{equation}
\begin{aligned}
\min_{\{u_k\}_{k=0}^{N_p-1}} & e_{N_p}^\top Q_N e_{N_p} + \sum_{k=0}^{N_p-1} \left( e_k^\top Q e_k + u_k^\top R u_k \right) , \\
\text{subject to} & \quad \quad e_k = r_k - C z_k , \\
& z_{k+1} = A z_k + B u_k, \quad k = 0, 1, \ldots, N_p - 1 , \\
& \quad z_{\min} \leq z_k \leq z_{\max} , \\
& \quad u_{\min} \leq u_k \leq u_{\max} , \\
& \quad \text{parameters} \quad z_0 = \Phi(x_0) , \\
& \quad r_k = \text{given} , \quad k = 0, 1, \ldots, N_p ,
\end{aligned}
\label{eq:mpc_cost}
\end{equation}
where \( u_{\min/\max} \) are the input bounds, \( z_{\min/\max} \) are the bounds on the lifted states, and \( r_k \) is the reference signal. Only the first control input \( u_0 \) from the optimal sequence is applied, and the optimization problem is resolved at the next time step.

\subsection{Iterative Data Augmentation}
Models trained on trajectories that are generated by randomly sampling the phase space and control input could perform poorly when subjected to \emph{optimal} feedback control input. 
This is mainly because these trained models lack exposure to the desired scenarios during the training process, especially for a relatively high-dimensional system.
To address this issue, we propose a straightforward iterative data augmentation technique inspired by the work of Uchida and Duraisamy~\cite{uchida2023control}. The idea is to iteratively collect the response of the true system under control input synthesized with the previous data-driven model.
Such response together with control input is used to augment the original dataset to further refine and improve the data-driven model. 
Our rationale is that such an experience from the true system can be used to better explore the phase space of the nonlinear system than merely ad hoc random sampling of the phase space, which is subject to the curse of dimensionality.

\begin{figure*}[!t]
    \centering
    \includegraphics[width=\linewidth]{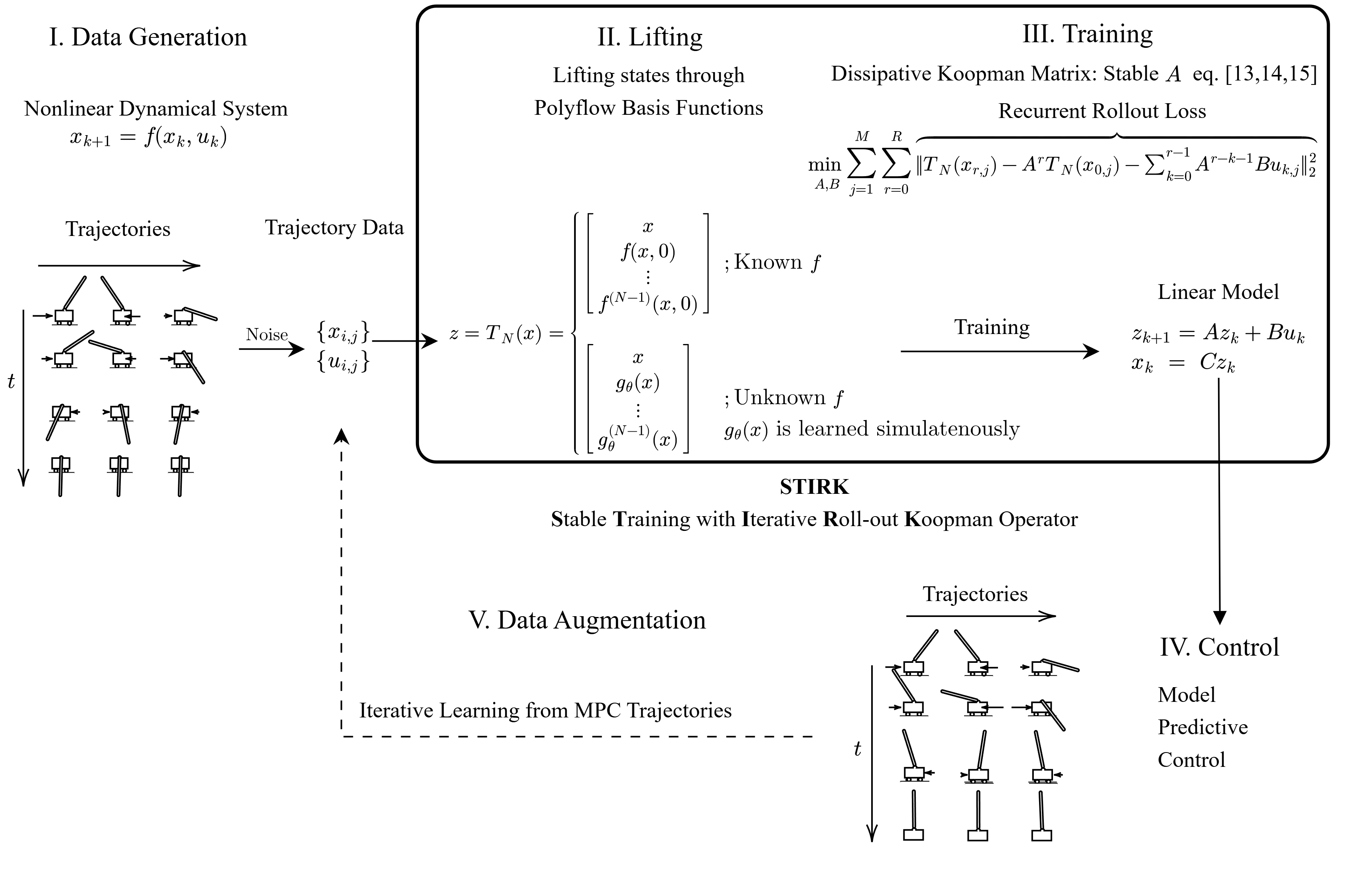}
    \caption{Overview of the \textbf{S}table \textbf{T}raining with \textbf{I}terative \textbf{R}oll-out \textbf{K}oopman Operator framework. The process begins with trajectory data collected from a nonlinear dynamical system. The states are lifted through Polyflow basis functions in \cref{eq:polyflow} if the system dynamics is known, or a recursive structure in \cref{eq:unknown_recursive} if the system dynamics is unknown, to obtain a higher-dimensional representation. This representation is then used in the roll-out EDMD algorithm to derive a linear model. The resulting linear model is employed for Model Predictive Control (MPC), which is further retrained through iterative data augmentation from closed loop trajectories.}
    \label{fig:summary}
\end{figure*}

Hence, our proposed framework is named as \textbf{STIRK} (\textbf{S}table \textbf{T}raining with \textbf{I}terative \textbf{R}oll-out \textbf{K}oopman Operator) as illustrated in figure \cref{fig:summary}.

\section{Numerical Experiments}
\label{sec:Numerical Examples}

In this section, we evaluate our proposed framework against other popular approaches in five different systems.

\begin{enumerate}
    \item \textbf{Autonomous Van der Pol Oscillator} – Used to compare prediction performance.
    \item \textbf{Cartpole System} – Evaluated on the task of stabilizing the pole in the upright position from \cite{underactuated}
    \item \textbf{Mountain Car (continuous action space version)} from \cite{towers2024gymnasium}
    \item \textbf{Lunar Lander (continuous action space version)} from \cite{towers2024gymnasium}
    \item \textbf{Panda Emika Robot} from \cite{gallouedec2021pandagym} 
\end{enumerate}

While the first two examples assume known system dynamics, the last three environments are evaluated under the assumption of \textit{unknown system dynamics}.


\subsection{Autonomous Van der Pol Oscillator}

\begin{figure}[!t]
    \centering
    \includegraphics[width=1\linewidth]{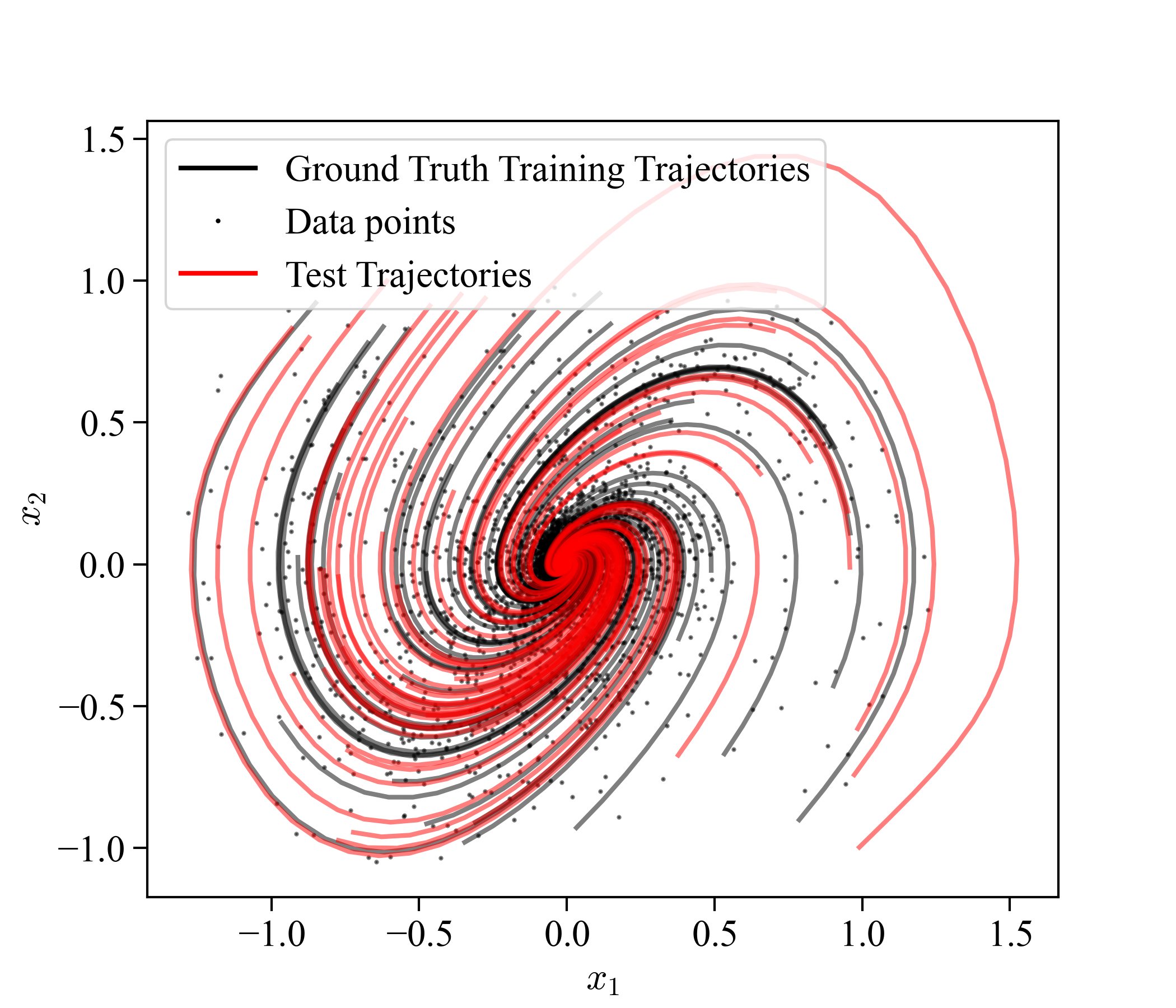}
    \caption{Visualization of the Van der Pol system's dynamical behavior, illustrating ground truth training trajectories, noisy data points (noise level 0.0599), and test trajectories in the $x_1,x_2$ plane.}
    \label{fig:data_distribution}
\end{figure}

We consider the autonomous Van der Pol Oscillator described by the equations:
\begin{align}
\dot{x}_{1} & = -x_{2}, \\
\dot{x}_{2} & = \mu \left(-1 + x_{1}^{2}\right) x_{2} + x_{1},
\end{align}
where \(\mu\) is set to 1.0. 
The Van der Pol Oscillator is a nonlinear system characterized by a non-conservative force that can either dampen or amplify oscillations depending on the state.
Since our framework relies on the discrete-time dynamical system, the continuous-time version of the van der Pol system is discretized using a fourth-order Runge-Kutta method with a time step \(\Delta t\) of 0.1. 
The training dataset is generated by running simulations for 100 steps, also with a time step of \(\Delta t\) of 0.1. 
After simulation, Gaussian noise is added to the states, with 10 logarithmically spaced noise levels from \(10^{-3}\) to \(10^{-1}\). 
For models trained on a single trajectory, the initial condition is sampled from a normal distribution with a mean of \([-1.0, -1.0]\) and a standard deviation of \([0.05, 0.05]\). 
For models trained on multiple trajectories, 50 initial conditions are uniformly sampled from the range \([-1, 1]\). To ensure robustness and generalizability, we created 10 sets of training data using different random seeds. Figure \cref{fig:data_distribution} illustrates one set of training trajectories and test trajectories within the $x_1,x_2$ plane.

\paragraph*{Prediction Performance Evaluation} 
We compare the proposed framework with the following noise-robust variants of DMD:
\begin{enumerate}
\item DMD~\cite{tu2013dynamic}, 
\item eDMD ~\cite{williams2015data} with Polyflow basis function and Radial Basis function,
\item OptDMD~\cite{askham2018variable} with Polyflow basis function, Radial Basis functions and Identity Observables, 
\item TLSDMD~\cite{ohmichi2022stable} with Polyflow basis function, Radial Basis functions and Identity Observables, 
\item FBDMD~\cite{dawson2016characterizing} with Polyflow basis function, Radial Basis functions and Identity Observables.
\end{enumerate}

To evaluate the prediction performance of these different algorithms, we use a normalized error metric that quantifies the discrepancy between the predicted trajectories and the ground truth trajectories. 
Let \( \mathbf{x}(t) \in \mathbb{R}^n \) represent the true state of the system at time \( t \) and \( \hat{\mathbf{x}}(t) \in \mathbb{R}^n \) represent the predicted state at time \( t \). The normalized prediction error over the trajectory is computed using the following formula:
\begin{equation}
\text{Error} = \frac{\|\hat{\mathbf{X}} - \mathbf{X}\|_F}{\|\mathbf{X}\|_F},
\end{equation}
where \( \|\cdot\|_F \) denotes the Frobenius norm, and \( \mathbf{X} \in \mathbb{R}^{T \times n} \) and \( \hat{\mathbf{X}} \in \mathbb{R}^{T \times n} \) are the matrices containing the true and predicted states over all time steps \( T \), respectively,
\begin{equation}
\mathbf{X} = \begin{bmatrix}
\mathbf{x_1} \\
\mathbf{x_2} \\
\vdots \\
\mathbf{x_T}
\end{bmatrix}, \quad \hat{\mathbf{X}} = \begin{bmatrix}
\hat{\mathbf{x}}_1 \\
\hat{\mathbf{x}}_2 \\
\vdots \\
\hat{\mathbf{x}}_T
\end{bmatrix}.
\end{equation}

\paragraph*{Training}
The Van der Pol system was initially trained using the Adam optimizer with a base learning rate of 0.01. A triangular-mode cyclic learning rate scheduler was employed, where the initial learning rate was halved each cycle. The maximum learning rate was set to 0.1, and each cycle increased for 500 epochs and decreased for another 500 epochs. The Polyflow order was configured to 4. The training process spanned 15,000 epochs with a batch size of 1,000. The rollout length was progressively increased, starting at 1 and doubling every 200 epochs, with a maximum rollout length initially set at 90. In the final 1,000 epochs, the optimizer was switched to L-BFGS with a learning rate of 0.01, applied on full batch data, to fine-tune the system.

\begin{figure*}[!t]
\centering
\subfloat[]{\includegraphics[width=3in]{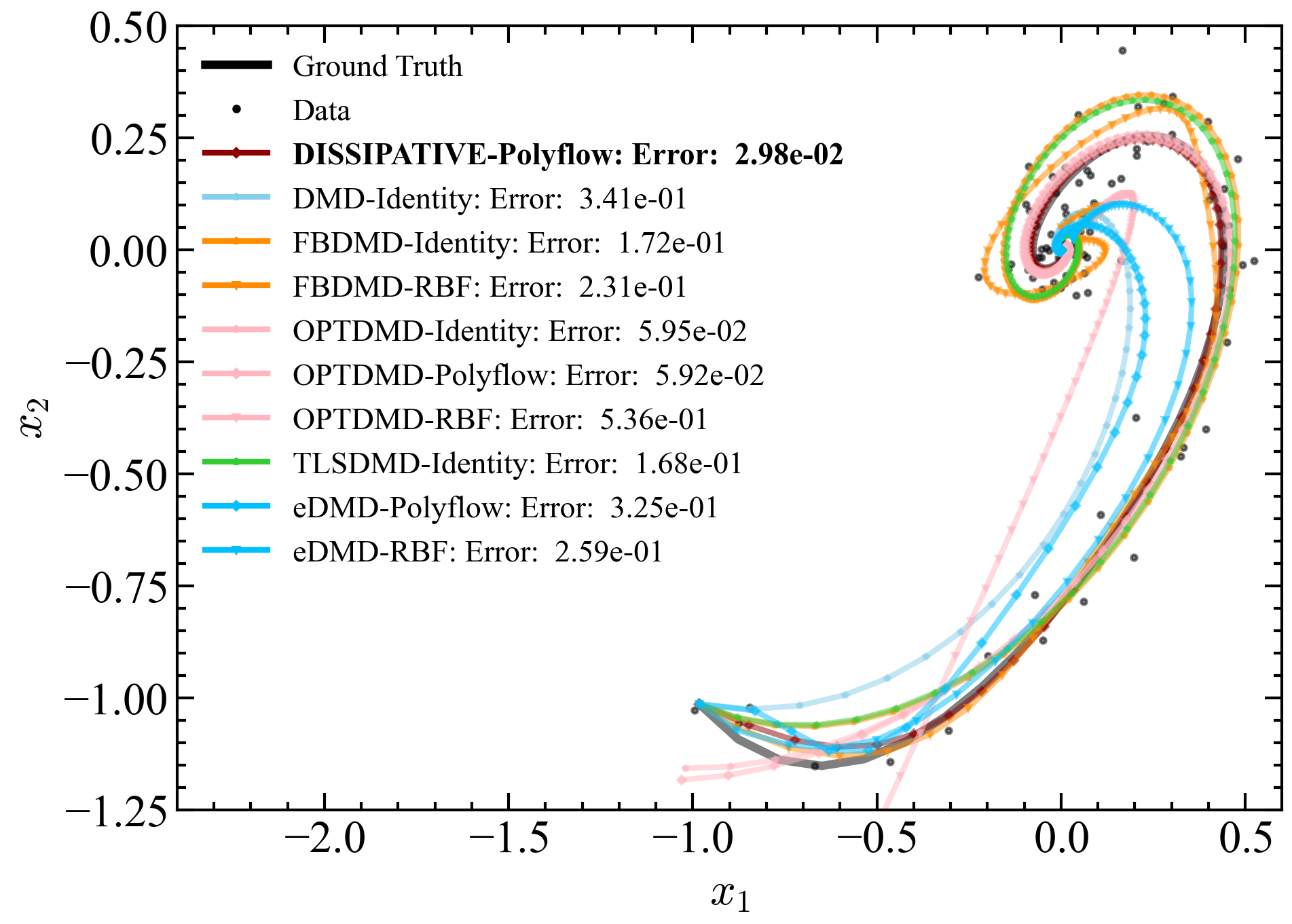}%
\label{fig:single_traj_prediction}}
\hfil
\subfloat[]{\includegraphics[width=3in]{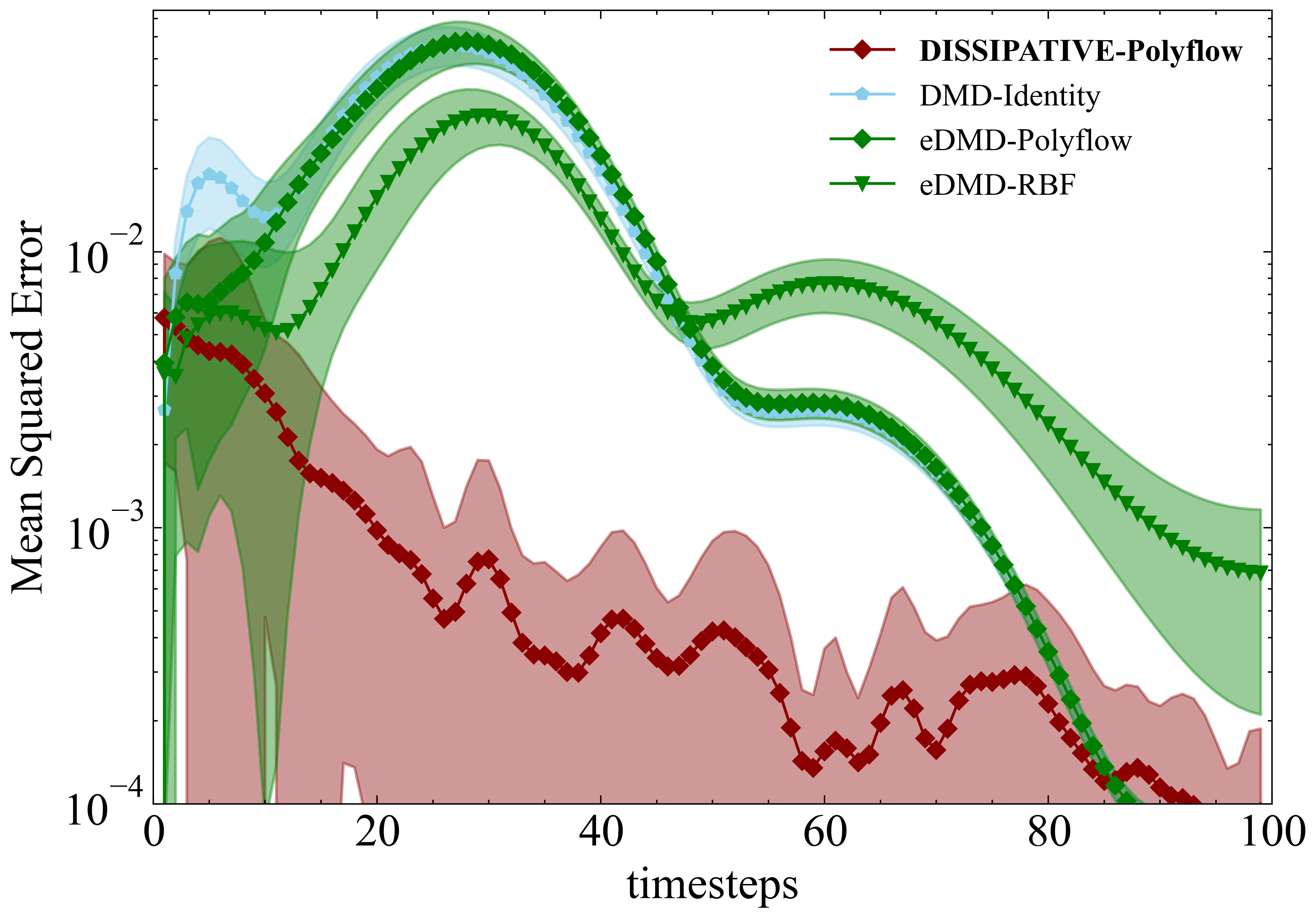}%
\label{fig:error_traj_sing_vdp_1}}
\caption{Comparison of model prediction performance for Van der Pol system. (a) Single Trajectory Prediction for Van der Pol System with Noise Level 0.0599. (b) Mean Squared Error (MSE) Trajectory over Timesteps for Van der Pol System with Noise Level 0.0599. The plot shows the evolution of MSE over 100 timesteps for different DMD-based methods. Shaded areas represent mean $\pm$ standard deviation of the mean squared errors. The Dissipative Koopman model with Polyflow basis functions maintains the lowest error throughout the timesteps, demonstrating robust performance under noisy conditions.}
\label{fig:single_traj_vdp_pred}
\end{figure*}

\paragraph*{Results for Single Trajectory}
Figure~\cref{fig:single_traj_prediction} compares the reconstruction performance of various methods based on dynamic mode decomposition (DMD) when trained on a single trajectory of the Van der Pol system at a noise level of 0.0599. The legend lists the normalized error for each method. The dissipative Koopman model with polyflow basis functions (Dissipative-Polyflow) achieved the lowest error of $2.98 \times 10^{-2}$, indicating a superior precision. OptDMD, using both identity and Polyflow basis functions, performed similarly well, closely approaching the performance of the dissipative Koopman model. However, since OptDMD also optimizes the initial condition, it shows a deviation from the actual initial condition. Other methods exhibited higher errors, indicating less accurate reconstructions under the given noise conditions.
Figure~\cref{fig:error_traj_sing_vdp_1} presents the mean squared error (MSE) trajectory over 100 time steps for the Van der Pol system at a noise level of 0.0599 for 10 different training sets. Each line represents the mean error for a different DMD-based method, with shaded areas indicating the standard deviation of the error.
The dissipative Koopman model with polyflow basis functions consistently maintains the lowest MSE across all timesteps, highlighting its robustness and accuracy under noisy conditions.
The eDMD with polyflow basis functions and the eDMD with RBF basis functions exhibit higher MSEs, particularly in the early timesteps, indicating less accurate predictions compared to the dissipative model.
The standard DMD with identity observables (DMD-Identity) also shows higher error compared to the dissipative Koopman model but performs better than some of the other methods.
The figures and \cref{tab:errors} demonstrate that the dissipative Koopman model with Polyflow basis functions outperforms other DMD-based methods in terms of prediction accuracy for the Van der Pol system under noisy conditions. This highlights the effectiveness of incorporating dissipative properties and structured observables in enhancing the robustness and accuracy of the model.

\begin{figure*}[!t]
\centering
\subfloat[Training Error]{\includegraphics[width=3in]{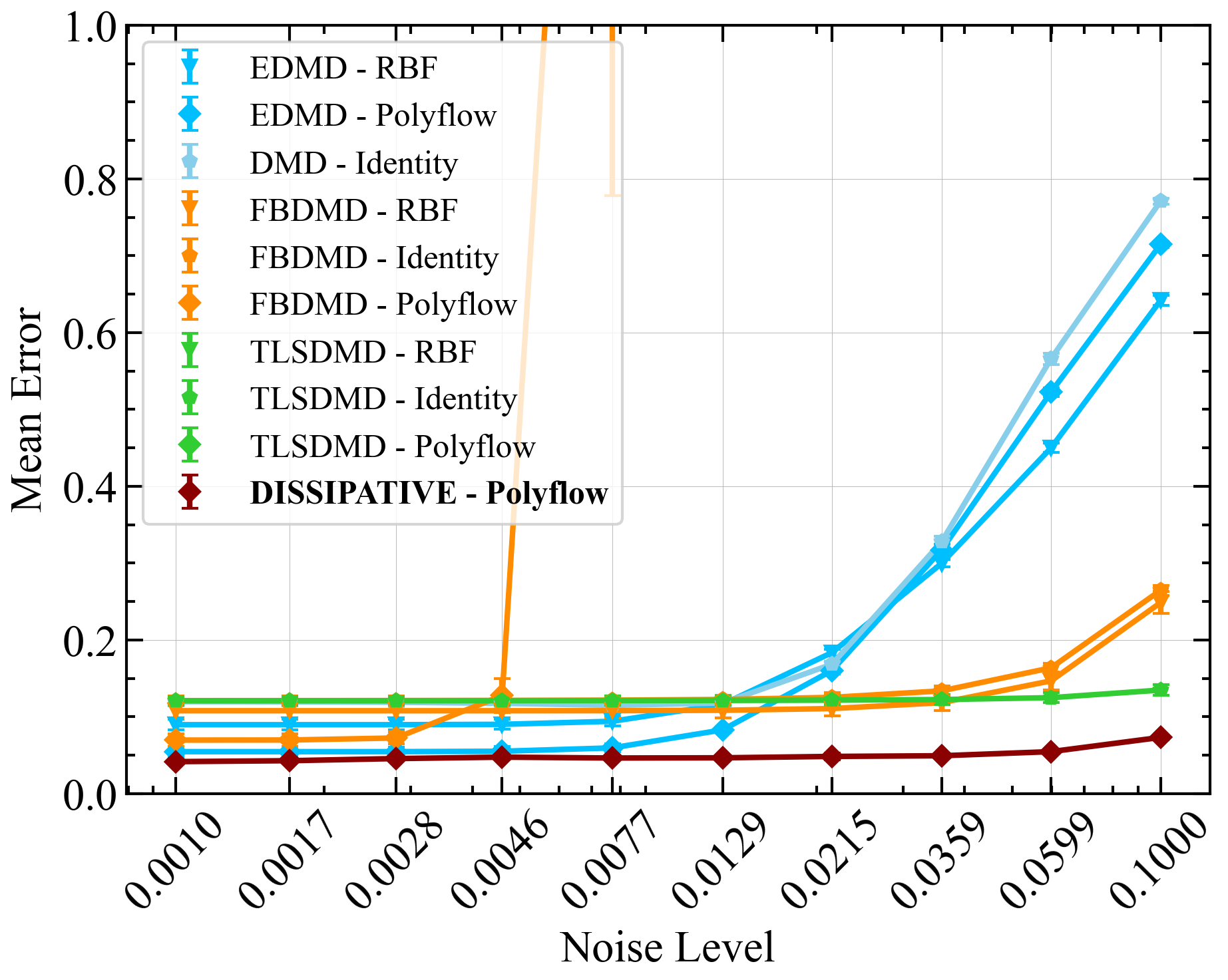}%
\label{fig:multiple_traj_vdp_comp}}
\hfil
\subfloat[Test Error]{\includegraphics[width=3in]{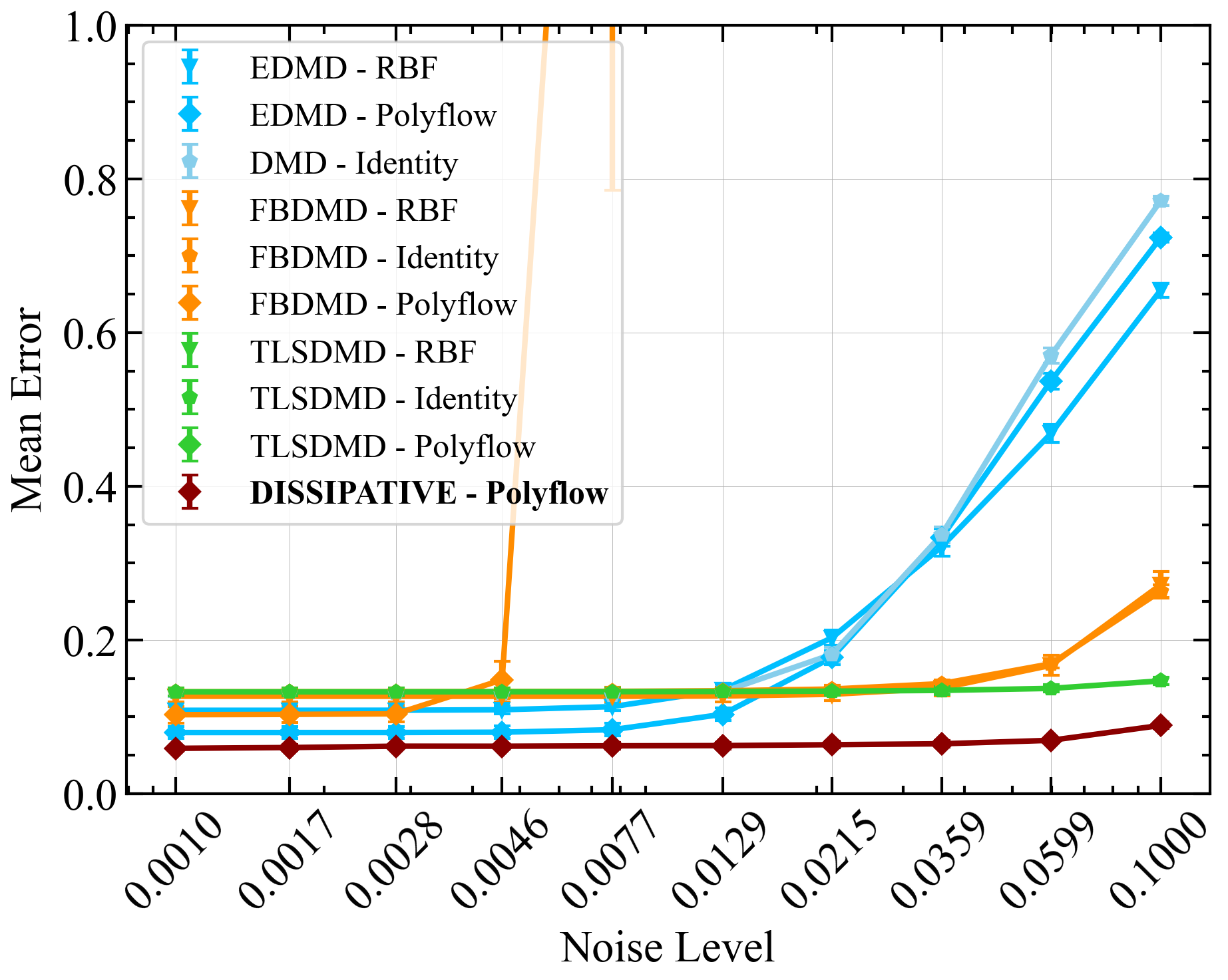}%
\label{fig:multiple_traj_vdp_comp_test}}
\caption{Multiple Trajectory Error Comparison for the Van der Pol (VDP) system. (a) Training Error. (b) Test Error.}
\label{fig:multiple_traj_vdp}
\end{figure*}

\paragraph*{Results for Multiple Trajectories}
Figure~\cref{fig:multiple_traj_vdp} compares the prediction errors for multiple trajectories of the Van der Pol (VDP) system across different noise levels. Panel (a) shows the training error, while panel (b) shows the test error. The dissipative Koopman model with polyflow basis functions consistently shows lower errors across various noise levels, demonstrating its robustness and accuracy in handling multiple trajectories.

\begin{table*}
\caption{Summary of Prediction Performance for Van der Pol System across different strategies}
\setlength{\tabcolsep}{4pt}
\label{tab:errors}
\begin{tabular}{ll|r|rrr|rrr|rr|c}
\hline
 & Strategy & DMD & \multicolumn{3}{c|}{FBDMD} & \multicolumn{3}{c|}{TLSDMD} & \multicolumn{2}{c|}{eDMD} & Dissipative (Ours) \\
 & \makecell{Lifting \\Function} & Identity & Identity & Polyflow & RBF & Identity & Polyflow & RBF & Polyflow & RBF & Polyflow \\
Error Type & Noise &  &  &  &  &  &  &  &  &  &  \\
\hline
\multirowcell{3}{Mean Normalized \\ Training Error} & 0.0010 & 0.1198 & 0.1211 & 0.0697 & 0.1079 & 0.1210 & NaN & NaN & 0.0547 & 0.0898 & \textbf{0.0417} \\
 & 0.0599 & 0.5660 & 0.1638 & NaN & 0.1468 & 0.1249 & NaN & NaN & 0.5225 & 0.4500 & \textbf{0.0547} \\
 & 0.1000 & 0.7711 & 0.2647 & NaN & 0.2486 & 0.1348 & NaN & NaN & 0.7149 & 0.6420 & \textbf{0.0733 }\\
\cline{1-12}
\multirowcell{3}{Mean Normalized \\ Test Error} & 0.0010 & 0.1315 & 0.1326 & 0.1028 & 0.1268 & 0.1325 & NaN & NaN & 0.0796 & 0.1085 & \textbf{0.0590} \\
 & 0.0599 & 0.5704 & 0.1691 & NaN & 0.1671 & 0.1370 & NaN & NaN & 0.5370 & 0.4688 & \textbf{0.0694} \\
 & 0.1000 & 0.7715 & 0.2634 & NaN & 0.2718 & 0.1468 & NaN & NaN & 0.7235 & 0.6552 & \textbf{0.0888} \\
\hline
\end{tabular}
\end{table*}

\subsection{CartPole System}
Let us consider the CartPole System \cite{underactuated},
\begin{align}
\dot{x}_1 &= x_3, \\
\dot{x}_2 &= x_4, \\
\dot{x}_3 &= \frac{u + m_p \sin x_2 \left(l x_4^2 - g \cos x_2\right)}{m_c + m_p \sin^2 x_2}, \\
\dot{x}_4 &= \frac{u \cos x_2 + m_p l x_4^2 \cos x_2 \sin x_2 - (m_c + m_p) g \sin x_2}{l (m_c + m_p \sin^2 x_2)}.
\end{align}

The CartPole system consists of a cart that moves along a track, balancing a pole attached to it via a pivot. The dynamics of the system are described by the equations above, where \(x_1\) represents the position of the cart, \(x_2\) represents the angle of the pole, \(x_3\) represents the velocity of the cart, \(x_4\) represents the angular velocity of the pole, and \(u\) is the control input force applied to the cart. The parameters of the system include the mass of the cart (\(m_c = 4kg\)), the mass of the pole (\(m_p = 1kg\)), the length of the pole (\(l = 1m\)), and the acceleration due to gravity (\(g = 9.81 m/s^2\)).
For our framework, the continuous-time CartPole system is discretized using a fourth-order Runge-Kutta method with a time step \(\Delta t\) of 0.05.
The training dataset is generated by running simulations for 200 steps and contains 50 different trajectories.
After the simulation, Gaussian noise is added to the states with noise level of 0.1. 
The initial conditions for the CartPole system are randomly sampled with the position (\(x_1\)) from \([-1, 1]\), the angle (\(x_2\)) from \([- \pi / 2, \pi / 2]\), the velocity (\(x_3\)) from \([-0.1, 0.1]\), and the angular velocity (\(x_4\)) from \([-0.1, 0.1]\).
Control inputs are generated using exponentially decaying sine waves of randomly sampled frequencies, amplitudes and phase angles. 
\paragraph*{Training}
The CartPole system was trained using the Adam optimizer with a learning rate of 0.001. The Polyflow order was set to 4. The training dataset consisted of 45 trajectories, while 5 trajectories were reserved for validation. The training process spanned 5,000 epochs with a batch size of 2,000. The maximum roll-out length was set to 32 and the roll-out length was increased every 300 epochs by a factor of 2.
\paragraph*{MPC Configuration}
The trained Koopman model is used as a surrogate model for model predictive control, following the approach in \cite{korda2018linear}. The prediction horizon, \(N_P\), is set to 20 steps. Control input constraints are defined as \(u_{\text{min}} = -20\) and \(u_{\text{max}} = 20\), with no constraints on the system states. For the cost function described in \cref{eq:mpc_cost}, \(Q_{N_P}\) and \(Q_{N}\) are set to ones on the first \(n\) diagonal elements, with zeros elsewhere. The parameter \(R\) is set to 0.

\begin{figure}[!t]
    \centering
    \includegraphics[width=1\linewidth]{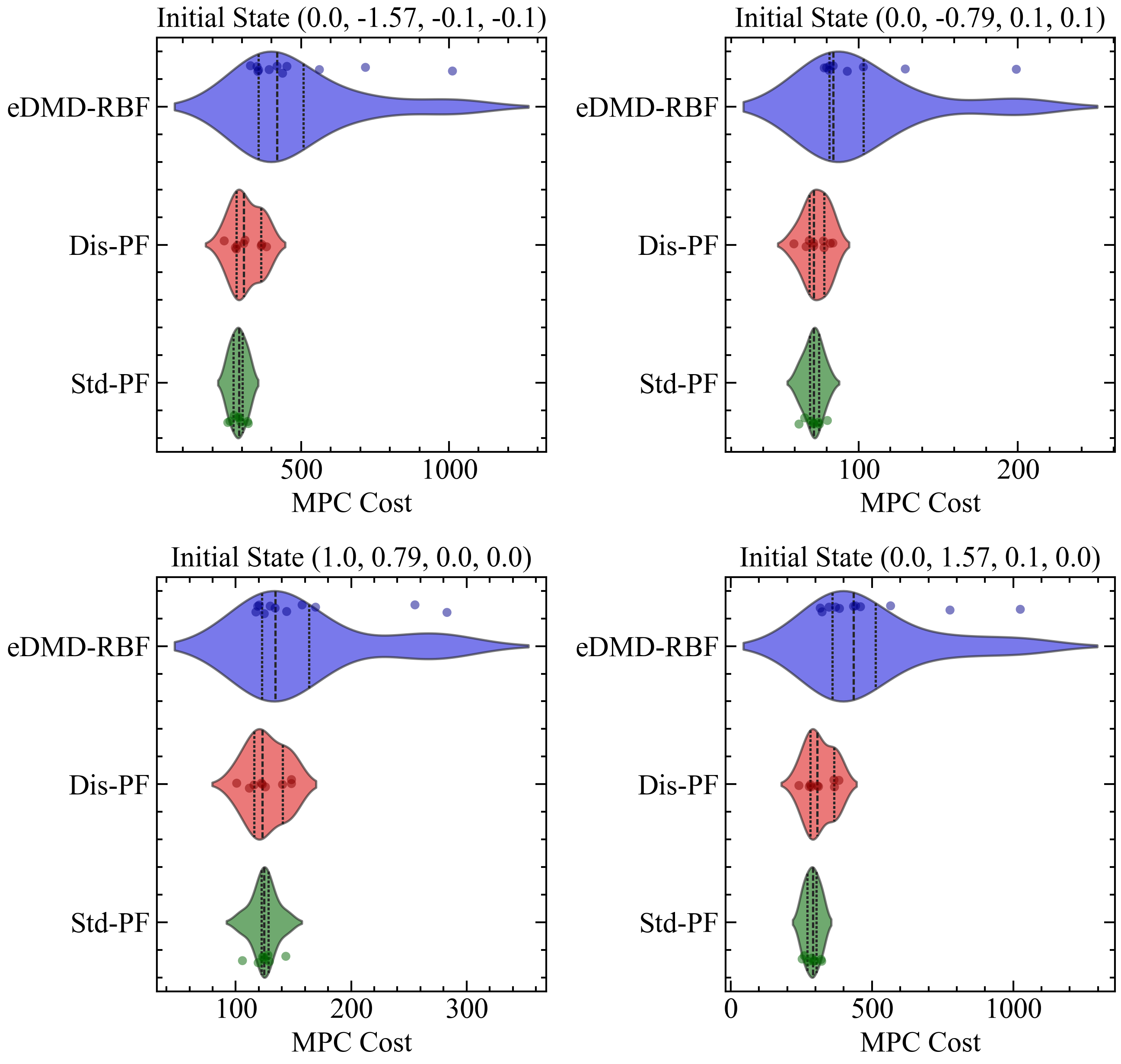}
    \caption{Violin plots showing the distribution of Model Predictive Control (MPC) costs for different models (eDMD-RBF, Dis-PF, Std-PF) under various initial states.}
    \label{fig:mpc_cost_vs_model}
\end{figure}

\paragraph*{Results}
Figure~\cref{fig:mpc_cost_vs_model} compares the Model Predictive Control (MPC) costs for the CartPole system using three different models under a noise level of 0.100. The models compared are EDMD with RBF basis functions (eDMD-RBF), the Dissipative Koopman model with Polyflow basis functions (Dis-PF), and the Standard Koopman model with Polyflow basis functions (Std-PF). The violin plots illustrate the distribution of MPC costs across various initial states.
The results indicate that both the Dissipative Koopman model with Polyflow basis functions (Dis-PF) and the Standard Koopman model with Polyflow basis functions (Std-PF) achieve relatively similar MPC costs. Both models show stable control performance across different initial states, demonstrating robustness in control.
In contrast, the EDMD with RBF basis functions (eDMD-RBF) exhibits higher variability in MPC costs and generally incurs higher average MPC costs compared to the Dissipative and Standard Koopman models. This suggests that while the EDMD-RBF model provides some level of stability, it does so at the expense of higher control costs and less consistency.
Overall, the comparison reveals that both the Dissipative and Standard Koopman models with Polyflow basis functions offer superior and more stable control performance compared to the EDMD-RBF model under the noise level of 0.1. This highlights the effectiveness of the Polyflow basis functions in achieving robust control in MPC applications.

\subsection{Mountain Car}
Mountain Car Continuous Version 0 is an environment from the Gymnasium library \cite{towers2024gymnasium} where a car is tasked with reaching the top of a mountain. The control action is a force applied along the x-axis with bounds \([-1.0 , 2.0]\). The observed state consists of the car's position along the x-axis and its velocity. 

We apply our learning framework to this environment under the assumption that the governing equation is unknown, and we have access only to noisy measurements of state-action sequences. For data generation, we use uniformly distributed random actions within the action space bounds and uniformly distributed initial positions between \([-0.6, -0.4]\) with zero initial velocity. The episode terminates upon reaching the target or after 1000 steps, with each step representing 0.02s. We collect 200 episodes of varying step lengths. 

To preprocess the dataset, we normalize the states and actions by subtracting their mean values and dividing by their standard deviations. Additionally, we introduce Gaussian noise with mean \(0\) and variances \(\sigma^2 = \{0.01, 0.1, 0.5\}\) to create three different training sets. A separate test dataset is generated following the same procedure.

\begin{figure}[!t]
\centering
\subfloat[Mountain Car Test Prediction Error]{\includegraphics[width=1\linewidth]{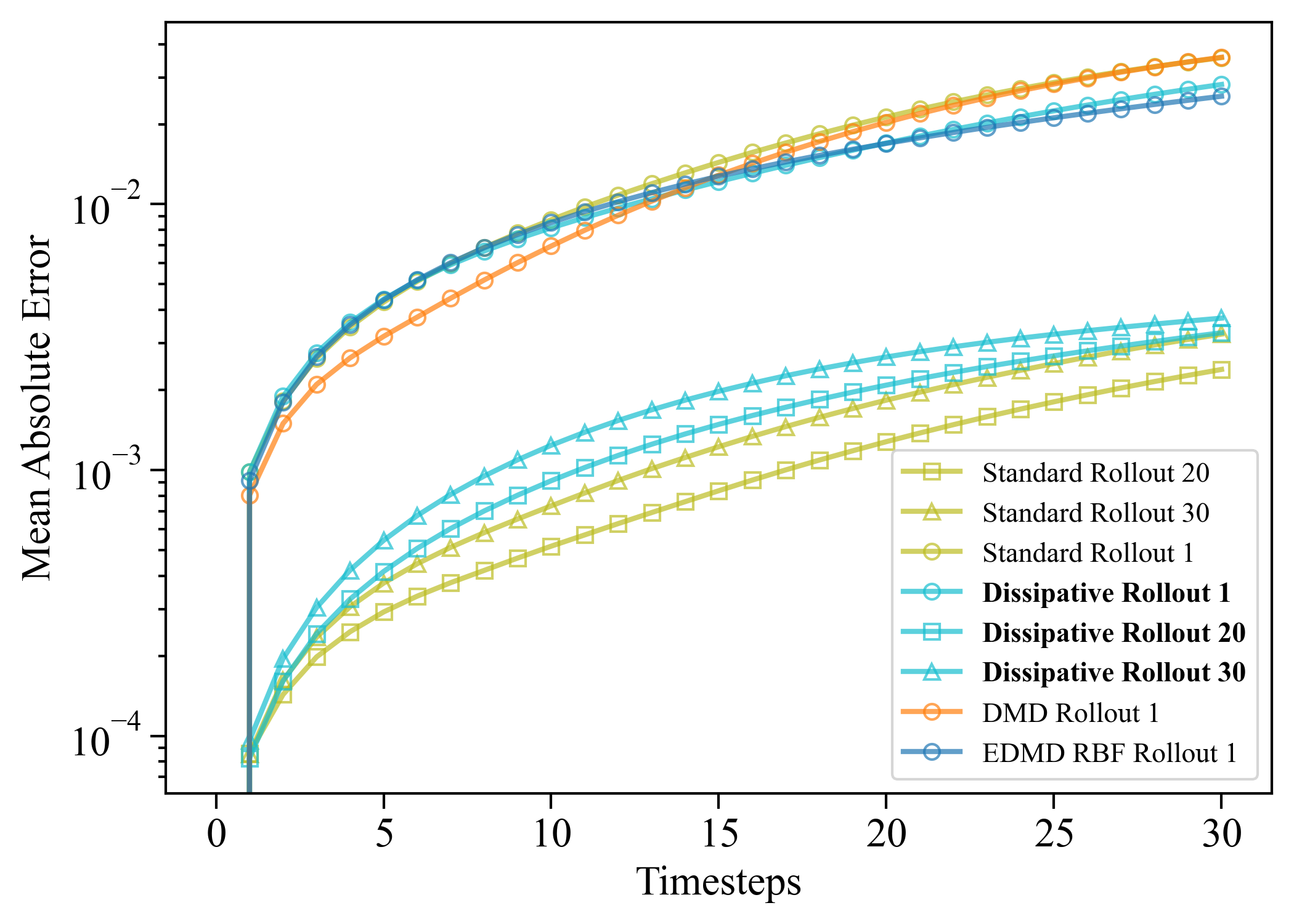}}
\vfill
\subfloat[Lunar Lander Test Prediction Error]{\includegraphics[width=1\linewidth]{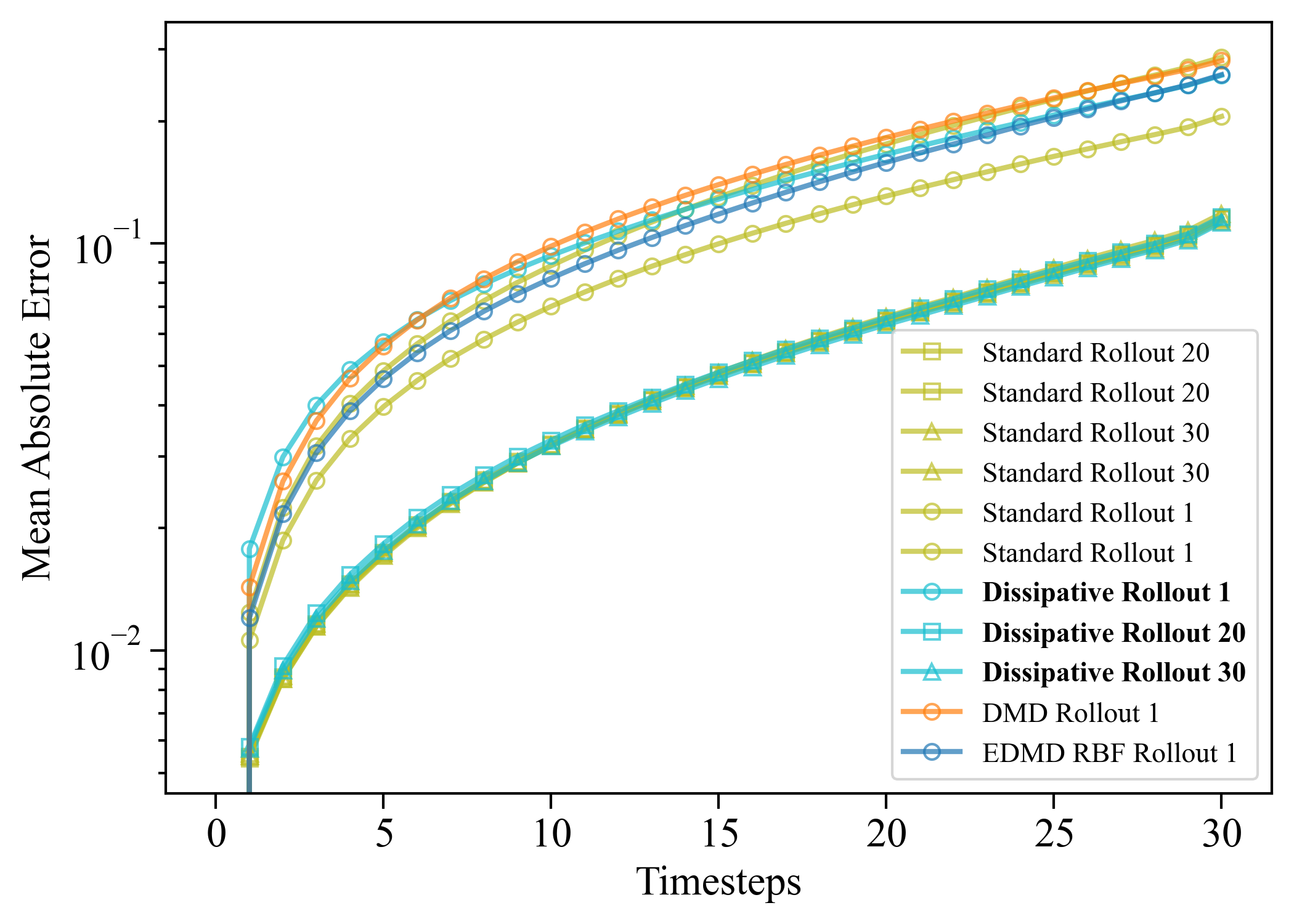}}
\vfill
\subfloat[Panda Test Prediction Error]{\includegraphics[width=1\linewidth]{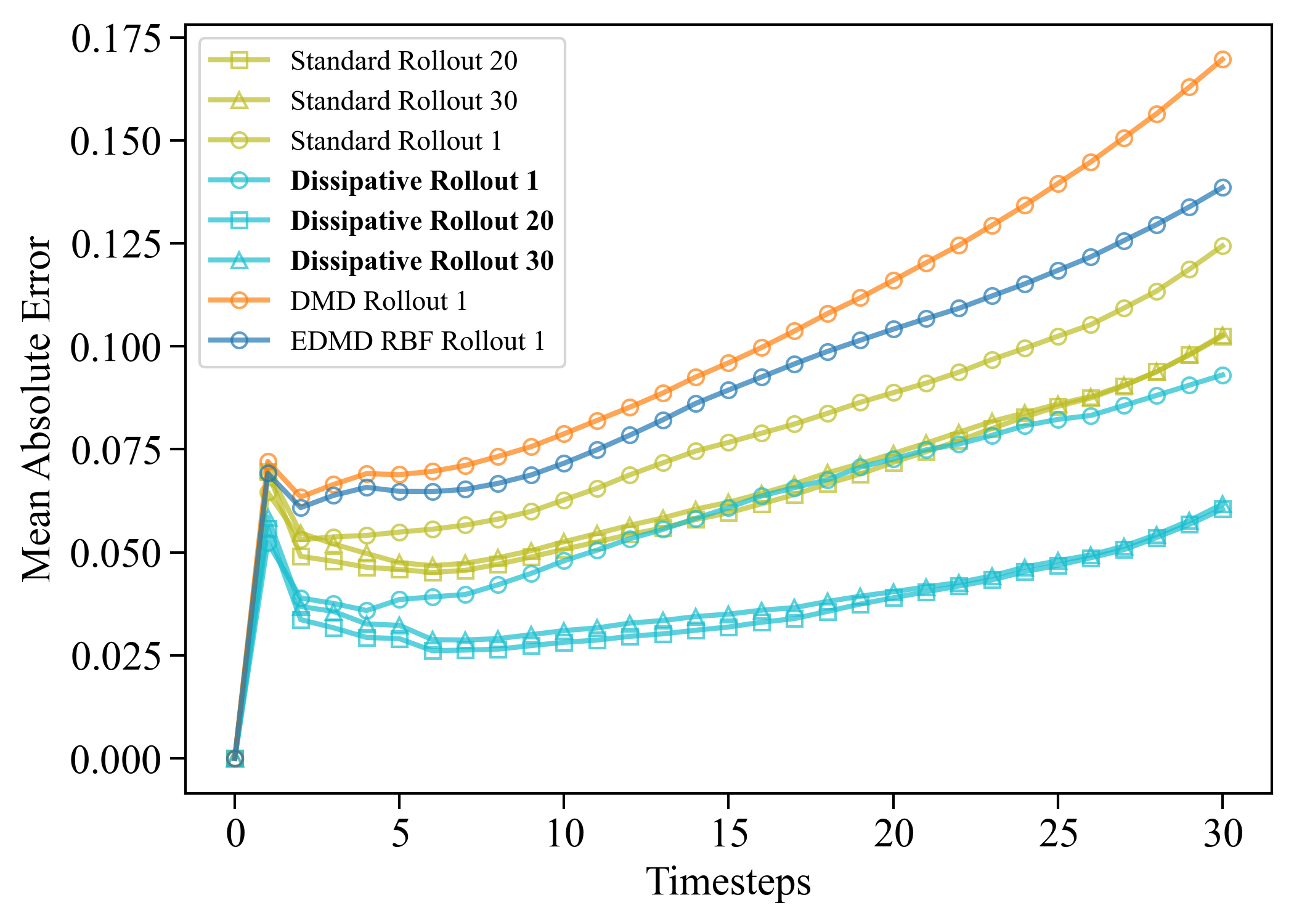}}
\caption{Test trajectory prediction errors for the Mountain Car, Lunar Lander, and Panda environments under unknown dynamics. For these experiments, a Polyflow order of 2 was used for the lifted state construction, and all models were trained with training data subject to 0.1 noise levels. Here, Adam optimizer was replaced by the SOAP optimizer \cite{vyas2024soap} with a learning rate of $3\times10^{-3}$.}
\label{fig:test-error}
\end{figure}

Test prediction errors are shown in \cref{fig:test-error}(a), MPC trajectories in \cref{fig:mpc-trajectory}(a), and noise sensitivity across sequence lengths in \cref{fig:heatmaps_noise} (left panel).

\subsection{Lunar Lander}
The Lunar Lander environment is a Gymnasium environment \cite{towers2024gymnasium} in which an agent controls the thrusters of a lander to achieve a soft landing on a designated pad. The state space originally comprises the position, velocity, angle, and angular velocity of the lander, while the action space consists of continuous thrust controls.

\begin{figure}[!t]
\centering
\subfloat[Mountain Car MPC Trajectory]{\includegraphics[width=1\linewidth]{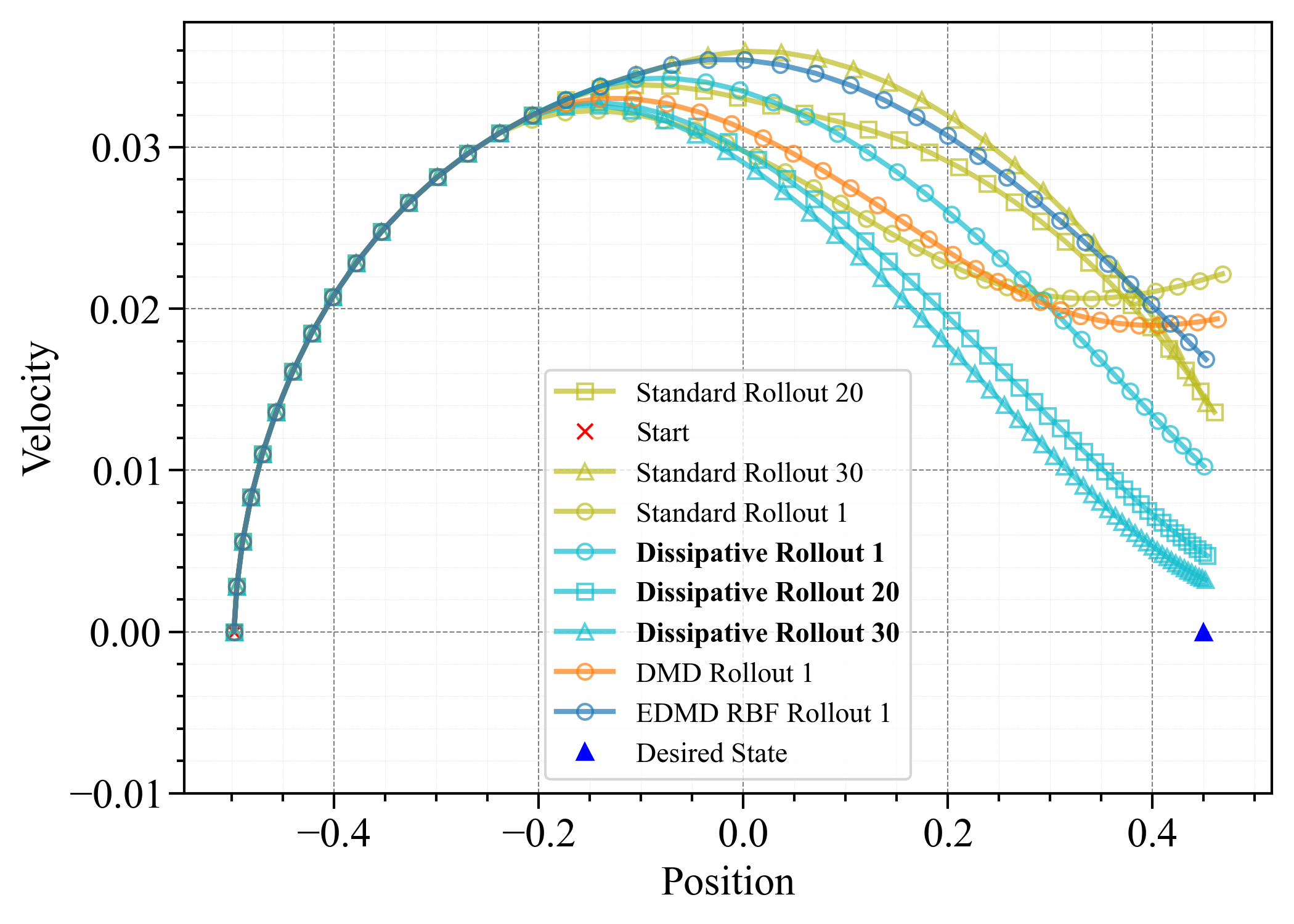}}
\vfill
\subfloat[Lunar Lander MPC Trajectory]{\includegraphics[width=1\linewidth]{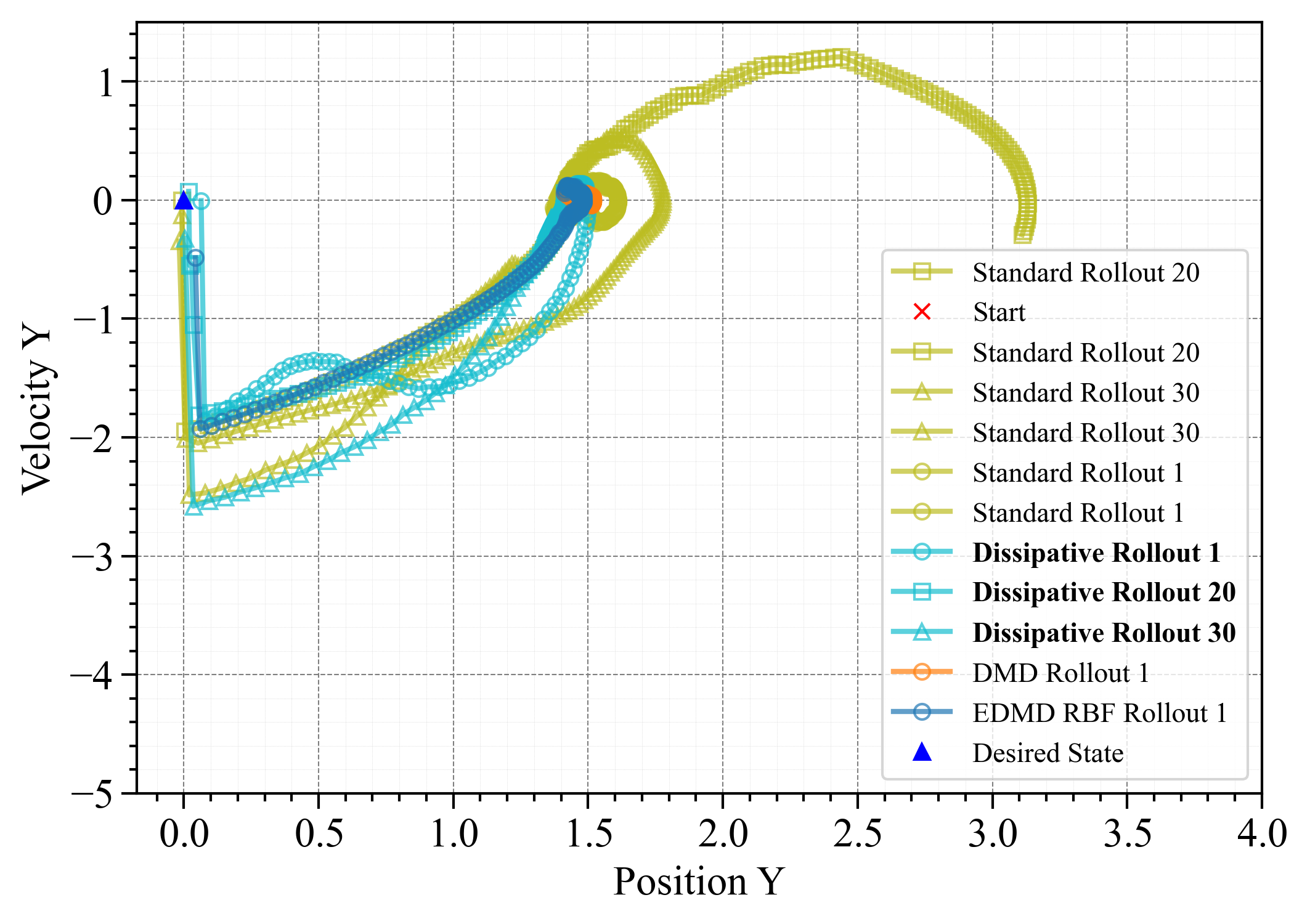}}
\vfill
\subfloat[Panda MPC Trajectory]{\includegraphics[width=1\linewidth]{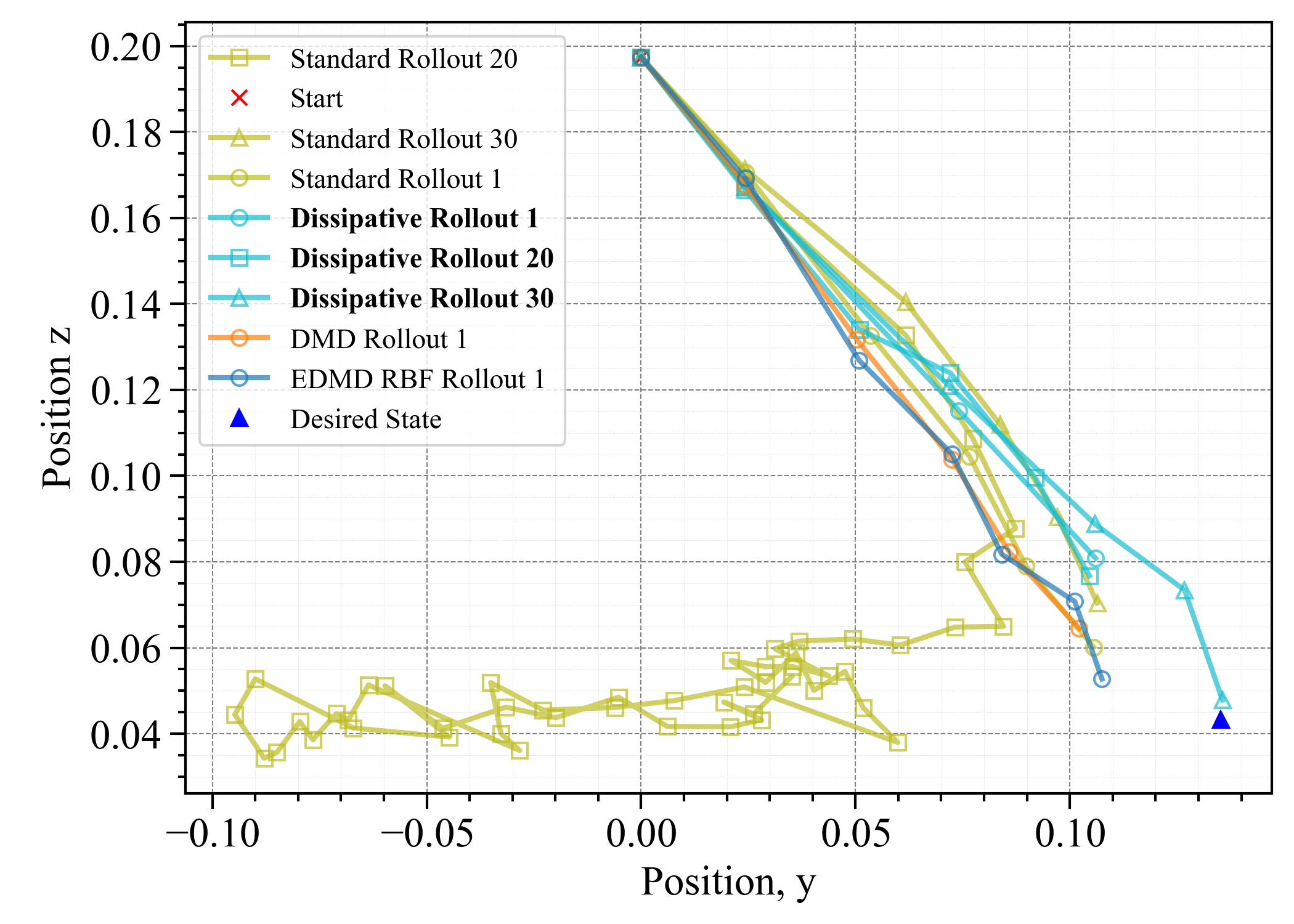}}
\caption{MPC-controlled trajectories for Mountain Car, Lunar Lander and Panda environments, trained on data with a noise level of 0.1. Model predictive control was implemented using the do-mpc \cite{FIEDLER2023105676} package.}
\label{fig:mpc-trajectory}
\end{figure}

The test prediction results are presented in \cref{fig:test-error}(b), and the resulting MPC-controlled trajectories are shown in \cref{fig:mpc-trajectory}(b). The heatmap in \cref{fig:heatmaps_noise} (middle panel) further illustrates the model’s robustness under different noise levels and rollout lengths.

\subsection{Panda Emika Franka Robot}
The Panda Emika Franka robotic arm, as provided in the panda-gym environment \cite{gallouedec2021pandagym}, is a 7-DoF manipulator widely used in industrial and research applications. In our experiments, the objective is to achieve precise target positioning by applying continuous control actions via a model predictive control (MPC) framework. The state space consists of the position and velocity of the robot's end effector, while the control inputs correspond to the forces applied at the end effector. 

\begin{figure*}[!t]
    \centering
    \includegraphics[width=1\linewidth]{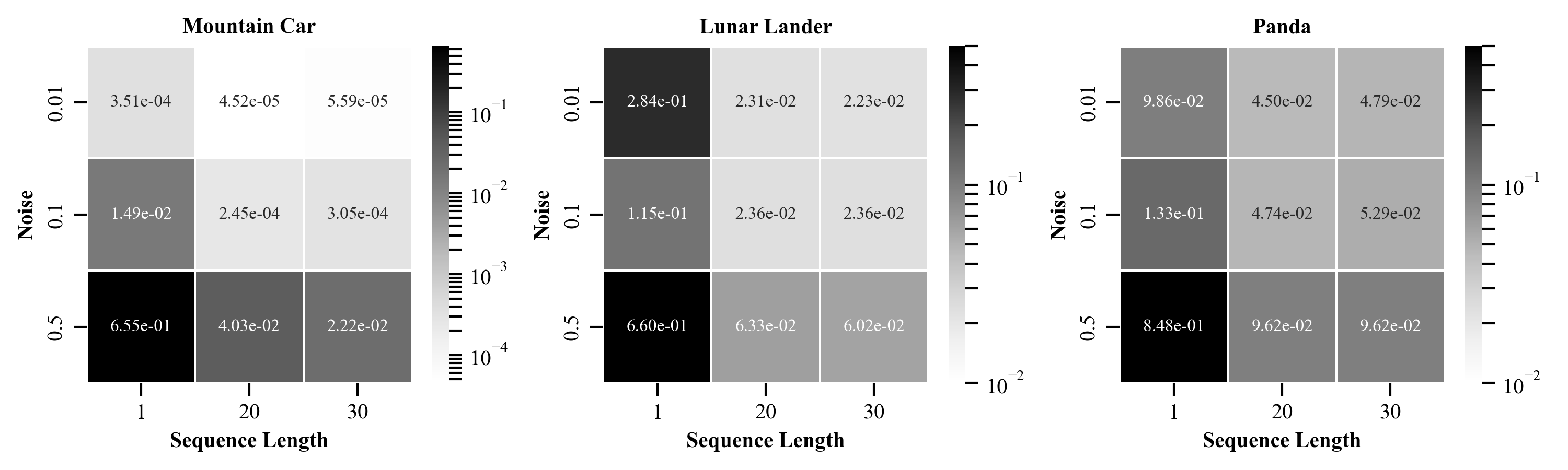}
    \caption{Heatmaps show predictive test performance across different noise levels and sequence lengths. Higher rollout length provides more robustness under noise in training data.}
    \label{fig:heatmaps_noise}
\end{figure*}

Prediction performance is visualized in \cref{fig:test-error}(c), and the corresponding MPC trajectories are shown in \cref{fig:mpc-trajectory}(c). Noise sensitivity under varying rollout lengths is depicted in the right panel of \cref{fig:heatmaps_noise}.

\subsection{Computational Cost Analysis with Varying Problem Dimensions}

To evaluate the computational scalability of our proposed Koopman learning framework, we benchmarked its performance on two nonlinear control environments: \textbf{Lunar Lander} (with a 6-dimensional state and 2-dimensional action space) and \textbf{Mountain Car} (with a 2-dimensional state and 1-dimensional action space). 

For each environment, we systematically varied:
\begin{itemize}
    \item the Polyflow order ($N_\phi$), and
    \item the rollout length ($R \in \{1, 10, 50\}$),
\end{itemize}
and recorded the training time (in milliseconds) and peak GPU memory usage (in megabytes) per epoch. All experiments were conducted using a single NVIDIA RTX 3080 GPU.

\begin{table}[h]
\centering
\caption{Training Time (ms) / GPU Memory Usage (MB) per Epoch on the \textbf{Mountain Car} Dataset}
\label{tab:cost_mc_full}
\begin{tabular}{c|ccc}
\toprule
\textbf{Polyflow Order} & \textbf{R = 1} & \textbf{R = 10} & \textbf{R = 50} \\
\midrule
2   & 14 / 26   & 95 / 102    & 719 / 443   \\
4   & 20 / 34   & 193 / 185   & 940 / 858   \\
10  & 25 / 59   & 337 / 435   & 1974 / 2106 \\
20  & 65 / 101  & 538 / 851   & 3504 / 4186 \\
\bottomrule
\end{tabular}
\end{table}

\begin{table}[h]
\centering
\caption{Training Time (ms) / GPU Memory Usage (MB) per Epoch on the \textbf{Lunar Lander} Dataset}
\label{tab:cost_mc_lunar}
\begin{tabular}{c|ccc}
\toprule
\textbf{Polyflow Order} & \textbf{R = 1} & \textbf{R = 10} & \textbf{R = 50} \\
\midrule
2   & 25 / 27   & 129 / 112    & 610 / 491   \\
4   & 31 / 36   & 233 / 202    & 835 / 939   \\
10  & 36 / 64   & 496 / 473    & 2455 / 2288 \\
20  & 96 / 110  & 820 / 925    & 4109 / 4538 \\
\bottomrule
\end{tabular}
\end{table}

As illustrated in \Cref{tab:cost_mc_full,tab:cost_mc_lunar}, both training time and memory usage increase with higher Polyflow orders and longer rollout horizons. While longer rollouts provide richer temporal context for capturing long-term system dynamics, they also lead to greater computational overhead.

\section{Iterative Data Augmentation}
\label{sec:adaptive data augmentation}
To explore the phase space of a nonlinear system more effectively, we introduce an iterative data augmentation technique. Instead of relying solely on ad hoc random sampling, we augment the training data with failure trajectories from closed-loop control and then retrain the Koopman operator. This approach addresses challenges highlighted by \cite{do2024practical}, where the authors emphasize the need to identify inputs that closely approximate the final closed-loop input and to use randomized signals that resemble it. However, finding suitable randomized control signals can be challenging and may not be adequate for control tasks. Therefore, we propose a systematic approach: augmenting the training data with control signals from previous failed control tasks in a closed-loop setting. This approach involves collecting the system's response to control inputs from data-driven based model predictive controller. The responses, along with the corresponding control inputs, are used to enhance the original data set, allowing further refinement and improvement of the model.
\begin{figure}[!t]
    \centering
    \includegraphics[width=1\linewidth]{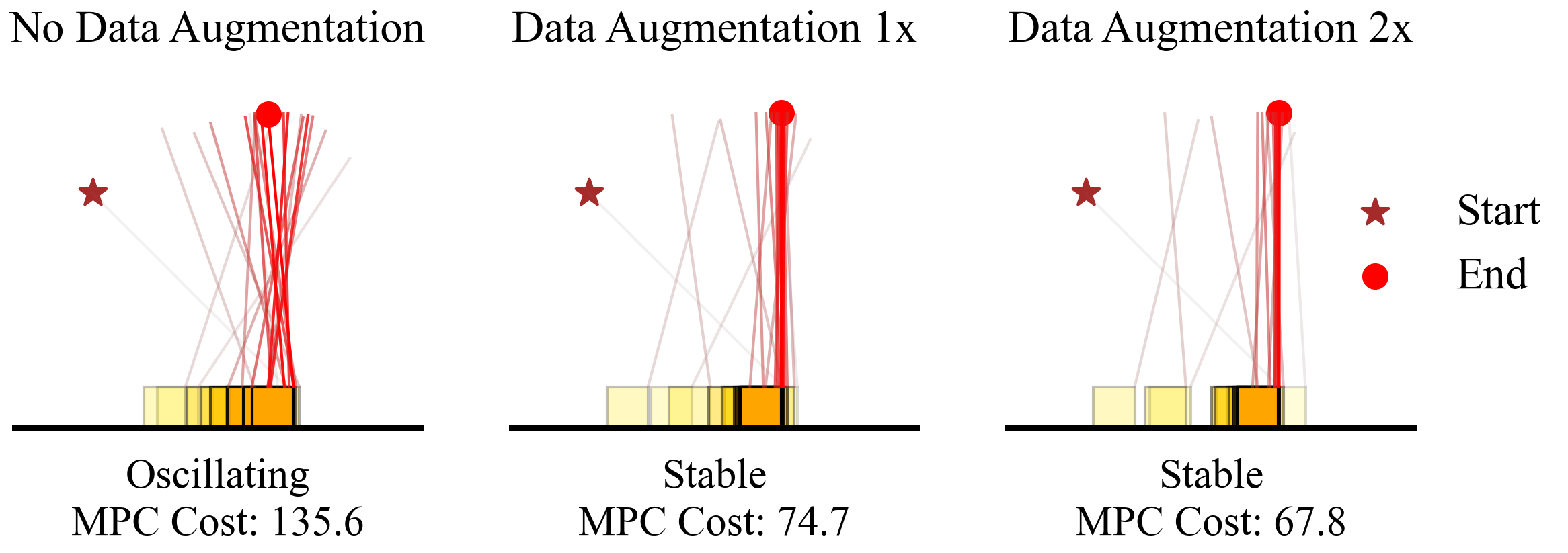}
    \caption{Trajectory visualization of the CartPole system over three learning iterations. The starting positions are indicated in light grey, while the final positions are highlighted in red.}
    \label{fig:adaptive_cartpole_anim}
\end{figure}
Figure~\cref{fig:adaptive_cartpole_anim} visualizes the trajectory of the CartPole system in three different scenarios: no data augmentation, one-time data augmentation, and twice data augmentation. In the first scenario, where no data augmentation was applied, the Koopman model fails to stabilize the CartPole after 200 steps. However, when the failed trajectory, along with other failed trajectories from different initial conditions, is added to the training data, the Koopman model demonstrates significantly improved control outcomes. Performance notably improves with data augmentation, transitioning from oscillating to stable trajectories. This enhancement is also reflected in the reduction of the Model Predictive Control (MPC) cost, with costs decreasing from 135.6 (no augmentation) to 74.7 (one-time augmentation) and 67.8 (twice augmentation).
\begin{figure}[!t]
    \centering
    \includegraphics[width=1\linewidth]{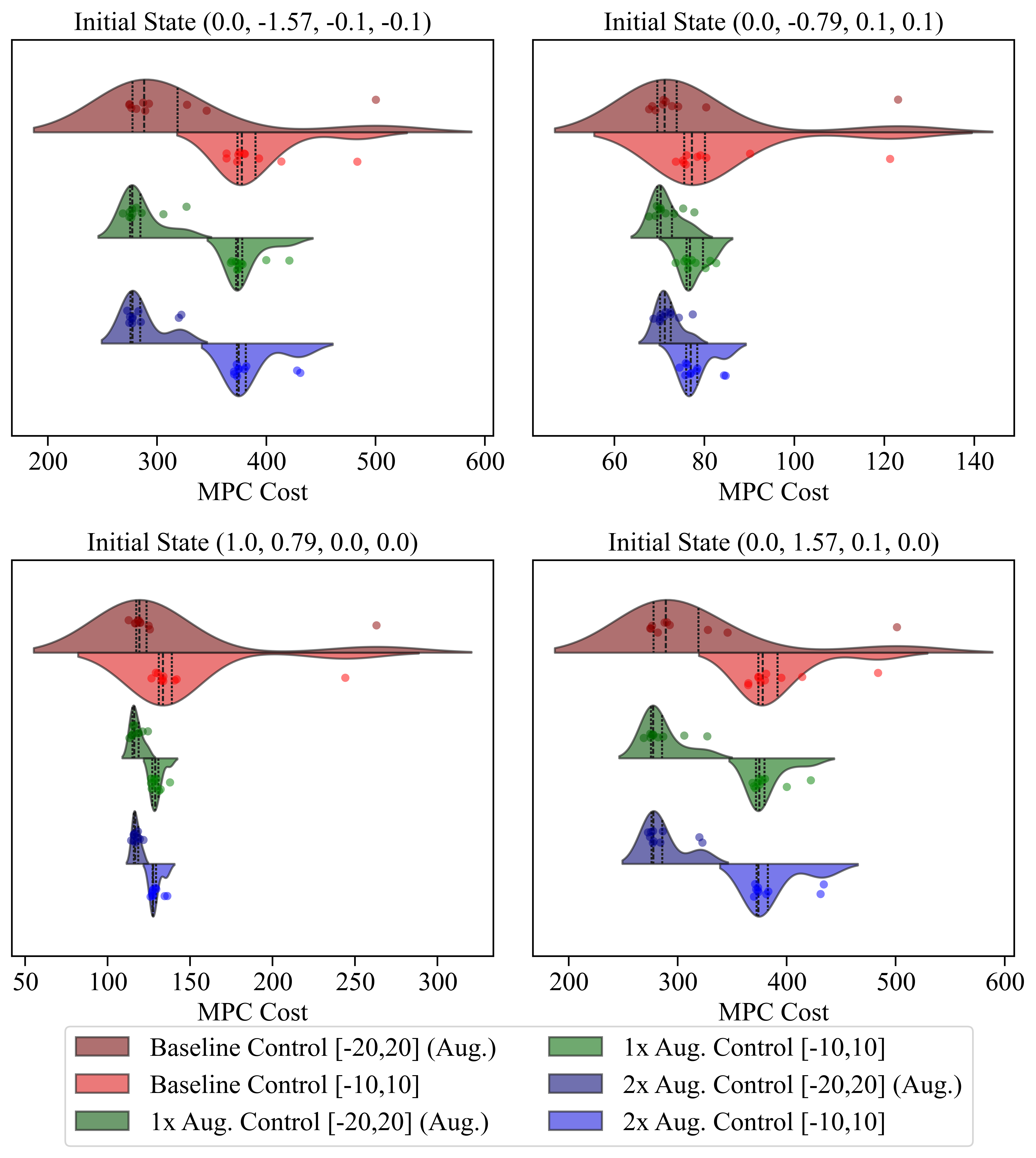}
    \caption{Violin plots depicting the distribution of Model Predictive Control (MPC) costs for CartPole problem for various tasks and iterative models under different initial states for noise level 0.1000.}
    \label{fig:adaptive_mpc_comparison}
\end{figure}
Figure~\cref{fig:adaptive_mpc_comparison} presents violin plots that illustrate the distribution of Model Predictive Control (MPC) costs for the CartPole system on two different tasks, seen and unseen, and iterative models under a noise level of 0.1. The ``seen" task refers to scenarios where input bounds are set between \([-20, 20]\), while the ``unseen" task refers to tighter input bounds of \([-10, 10]\). Data augmentation was performed using closed-loop trajectories from the seen task. It is important to note that iterative classical learning~\cite{mishra2010optimization} cannot be easily generalized to a different task. 
The violin plots display MPC costs for four different initial states and two levels of data augmentation, comparing the costs for both control tasks. Each subplot emphasizes the performance gains achieved through iterative data augmentation.
The results indicate that both seen and unseen control tasks benefit significantly from iterative data augmentation, with successive iterations leading to lower and more stable MPC costs. This shows the effectiveness of iterative data augmentation in improving the overall performance of the system.

\section{Ablation Studies}
\label{sec:ablation}
In this section, we will analyze the individual effect of the following components on the overall model performance with training data at different noise levels for the Van der Pol oscillator system.
\begin{itemize}
    \item Stabilized formulation of Koopman operator in Section \cref{sec:stabilized_formulation}, 
    \item Learning rate scheduler, 
    \item Optimizer switching,
    \item Progressive roll-out.
\end{itemize}
\begin{figure}[!t]
\centering
\subfloat[Effect of stabilized formulation of the Koopman operator.]{\includegraphics[width=1\linewidth]{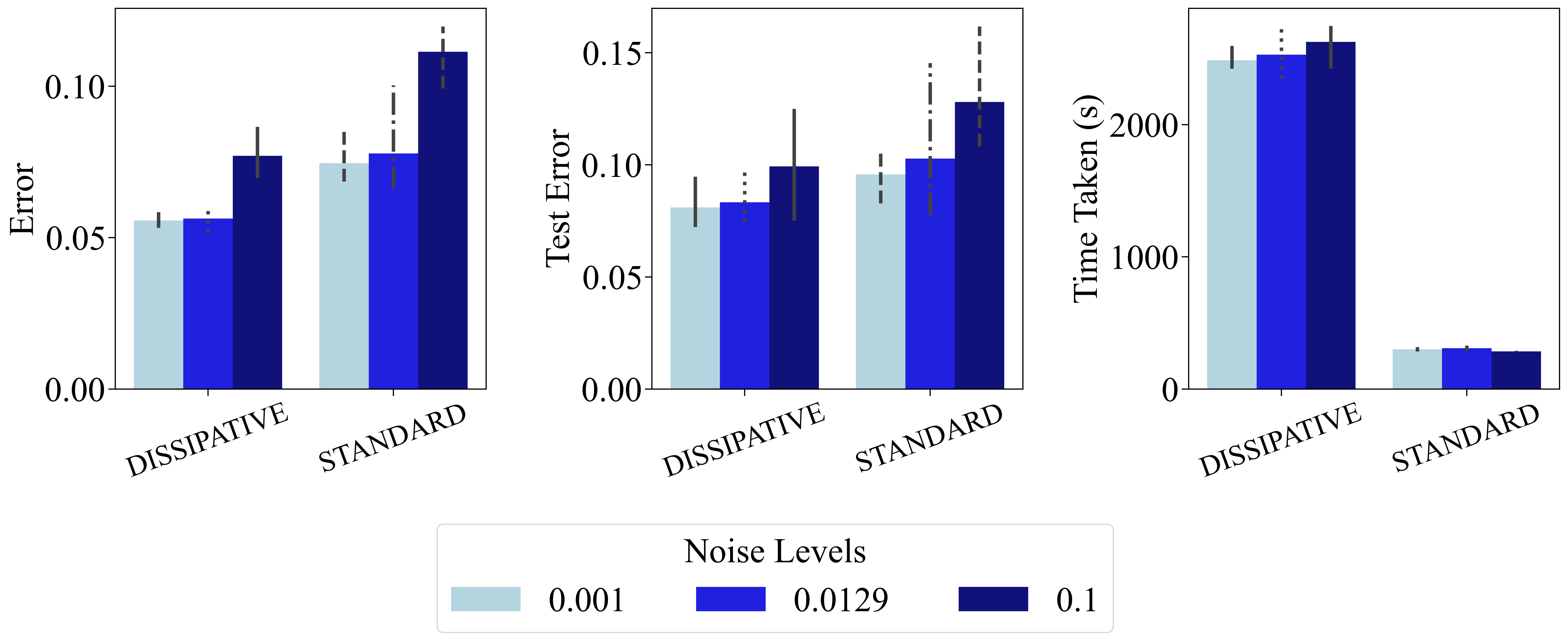}%
\label{fig:abl_1}}
\vfill
\subfloat[Effect of learning rate scheduler.]{\includegraphics[width=1\linewidth]{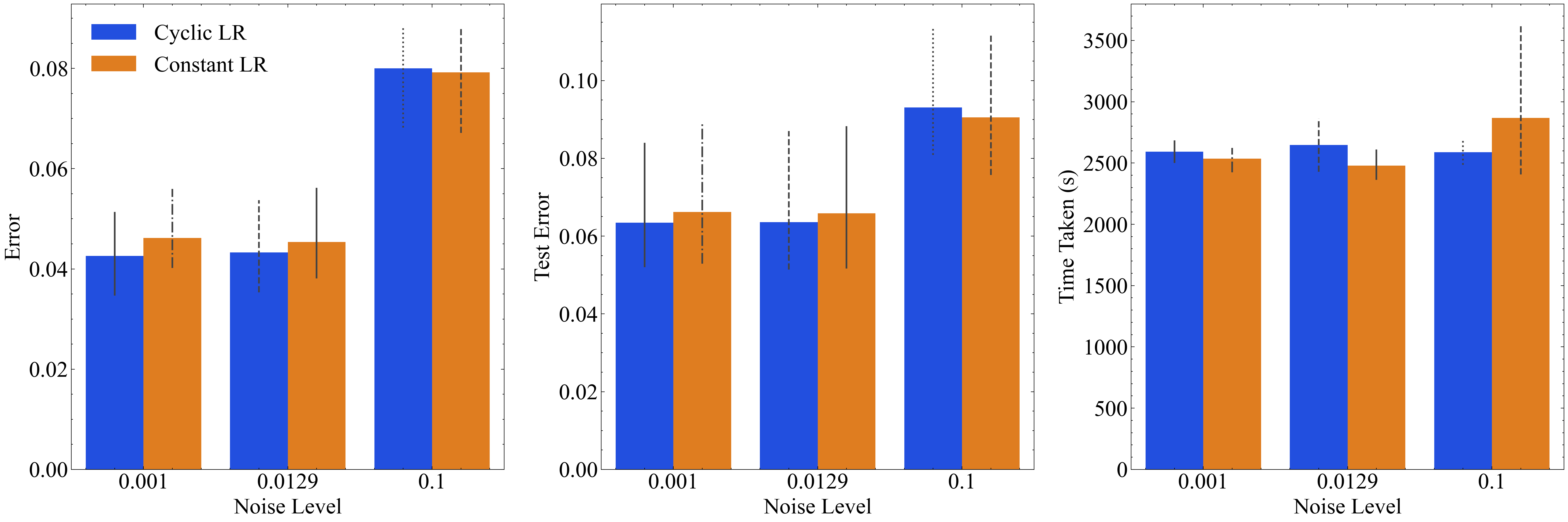}%
\label{fig:abl_2}}
\vfill
\subfloat[Effect of optimizer switching.]{\includegraphics[width=1\linewidth]{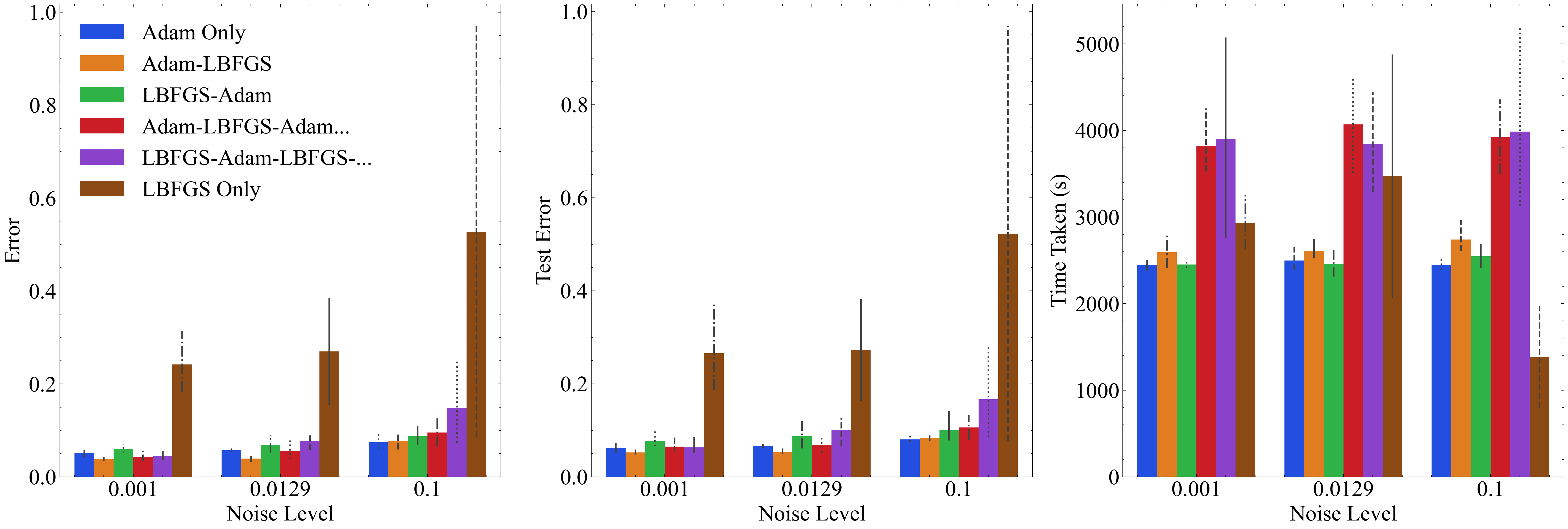}%
\label{fig:abl_3}}
\vfill
\subfloat[Effect of progressive roll-out.]{\includegraphics[width=1\linewidth]{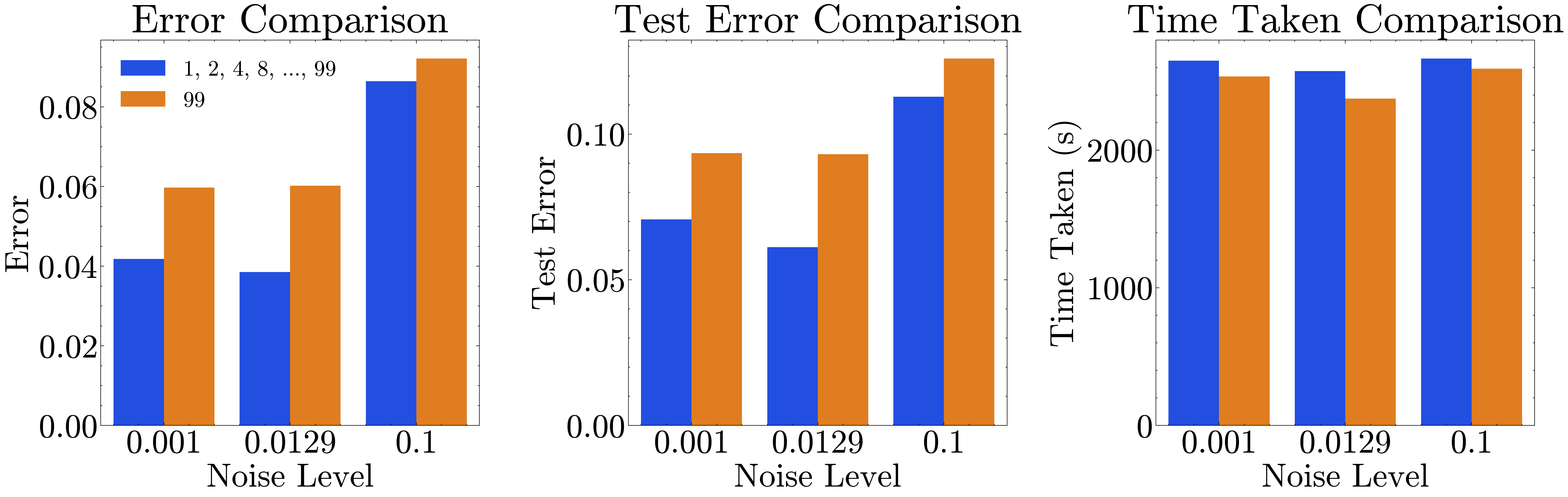}%
\label{fig:abl_4}}
\caption{Ablation study on the Van der Pol Oscillator system. (a) Effect of stabilized formulation of the Koopman operator. (b) Effect of learning rate scheduler. (c) Effect of optimizer switching. (d) Effect of progressive roll-out.}
\label{fig:abl}
\end{figure}

\subsection{Dissipative Koopman Model versus Standard Koopman Model}
\begin{figure}[!t]
    \centering
    \includegraphics[width=1\linewidth]{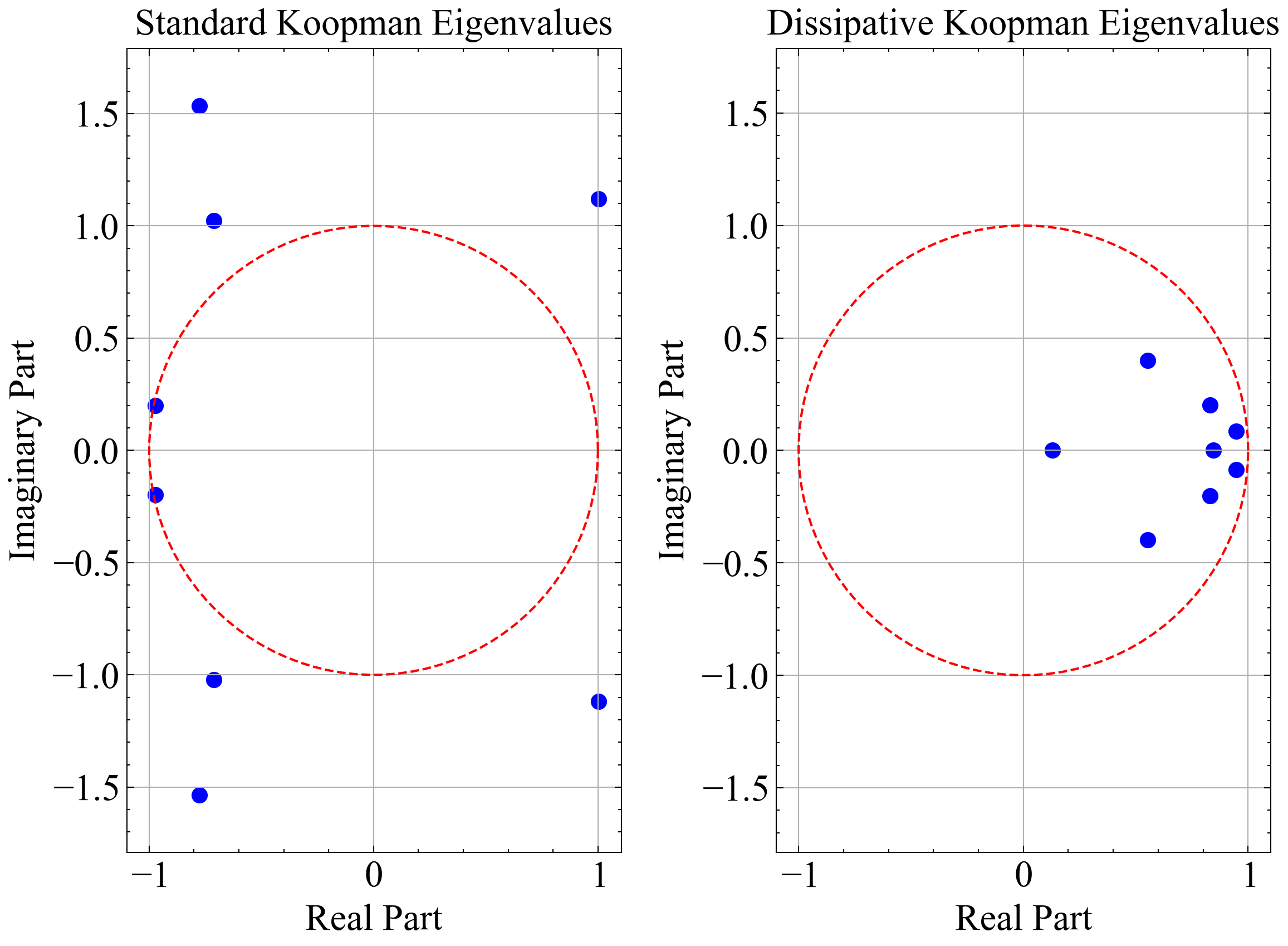}
    \caption{Eigenvalues of the Standard and Dissipative Koopman operators for the Van der Pol system with a noise level of 0.0599 after 1,000 training epochs}
    \label{fig:dis_vs_std}
\end{figure}

\begin{figure}[!t]
    \centering
    \includegraphics[width=1\linewidth]{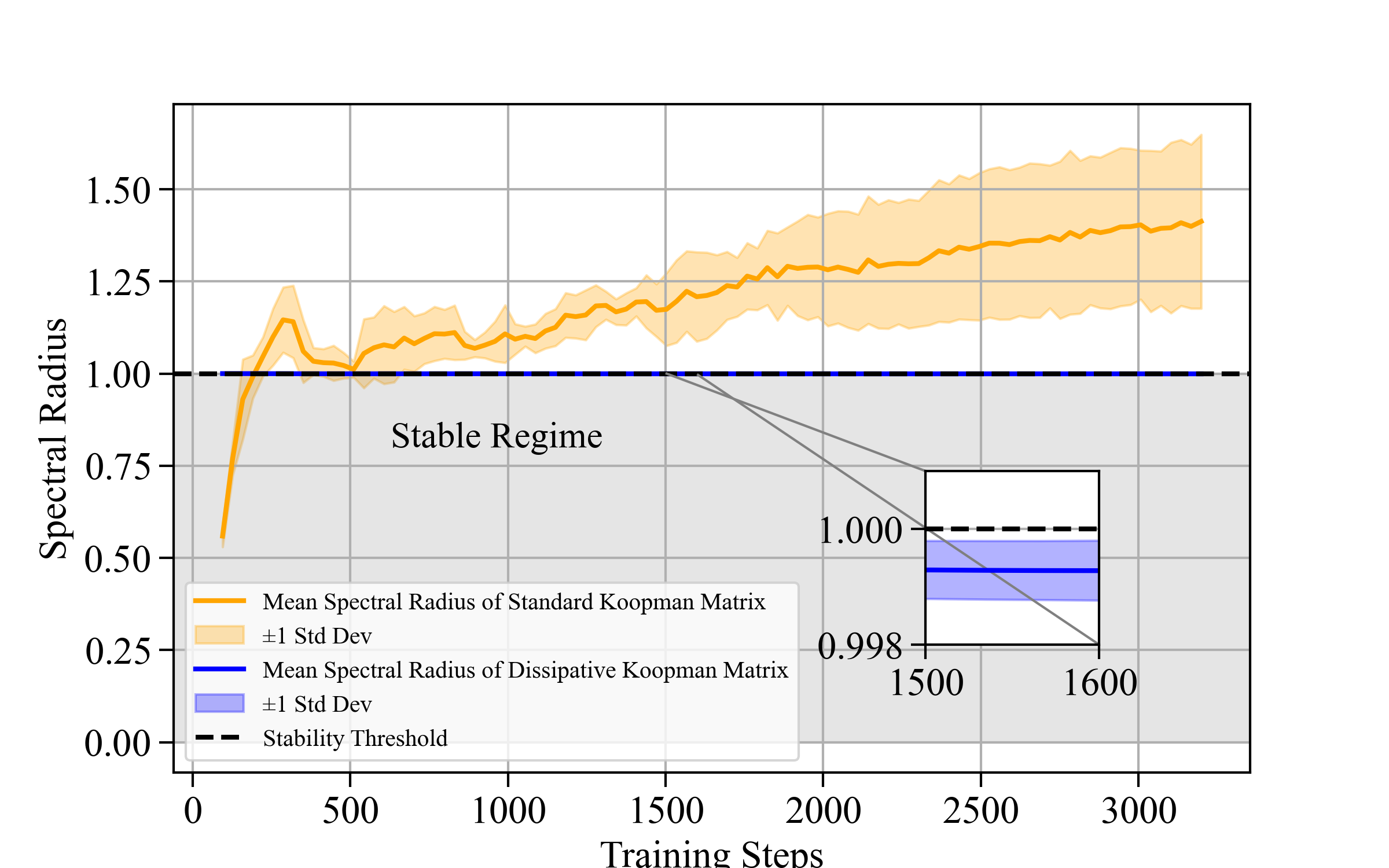}
    \caption{Evolution of the spectral radius for the Standard and Dissipative Koopman operators in the Panda environment with noise level 0.1, trained on sequences of length 20. The plot illustrates the stability behavior of both models over training steps. The dissipative Koopman operator consistently maintains its spectral radius below the stability threshold ($\rho = 1$), as highlighted in the inset, while the standard Koopman operator exceeds this threshold over time, indicating potential instability.}

    \label{fig:std_vs_dis}
\end{figure}

Figure~\cref{fig:abl_1} compares the performance of the dissipative Koopman model and the standard Koopman model in different noise levels for Van der Pol system. The Dissipative Koopman Model consistently shows lower mean errors and mean test errors, underscoring its superior prediction accuracy and robustness. However, the training time for the Dissipative Koopman Model is considerably longer than that of the Standard Koopman Model, due to the increased number of learning parameters and the computationally intensive operations, such as matrix exponentiation and pseudo-inverse computation.
Despite its relatively short training time, a significant drawback of the Standard Koopman Model is its sensitivity to the initial conditions of the Koopman matrix. If these initial conditions are not carefully tuned, the model can produce unstable eigenvalues, leading to poor performance. Figure~\cref{fig:dis_vs_std}  illustrates this issue by showing the discrete-time eigenvalues of the standard and dissipative Koopman operators for the Van der Pol Oscillator System at a noise level of 0.0599, using a different weight initialization than in previous cases. For fair comparison, both standard and dissipative Koopman models are initialized with stable Koopman matrices. However, after 1,000 training epochs, some of the discrete-time eigenvalues in the standard Koopman matrix move outside the unit circle, leading to a loss of stability, while the dissipative Koopman model never loses stability during training. Additionally, \Cref{fig:std_vs_dis} shows the evolution of the spectral radius ($\rho$, defined as the largest absolute value of the eigenvalues of the Koopman operator) over training steps for the Panda environment, reinforcing the long-term stability of the dissipative Koopman model compared to the standard one.

\subsection{Learning Rate Scheduler}
\Cref{fig:abl_2} compares the performance of the model using constant learning rate (LR) scheduling and cyclic learning rate scheduling under different noise levels for Van der Pol system. A \emph{cyclic learning rate} scheduler varies the learning rate between a predefined minimum and maximum in a periodic manner, often improving convergence stability \cite{smith2017cyclical}. For low noise, the cyclic learning rate scheduling achieves a 20\% reduction in training error and a 16.67\% reduction in test error, with a 9.09\% reduction in time taken compared to the constant learning rate. In medium noise scenarios, the cyclic learning rate scheduling shows a 16.67\% improvement in training error, 14.29\% in test error, and an 8.70\% improvement in time taken. However, at high noise levels, cyclic learning rate scheduling underperforms compared to the constant learning rate, with higher training and test errors, as well as slightly longer training times. These results indicate that while the cyclic learning rate scheduler is beneficial at low and medium noise levels, it is less effective in high-noise environments.

\subsection{Optimizer Switching}
Figure~\cref{fig:abl_3} examines the effect of switching optimizers on the model's performance for Van der Pol system. Different combinations of Adam and LBFGS optimizers were used, and their impact on mean errors, mean test errors, and training times was analyzed. Switching from the Adam optimizer to LBFGS after 4000 epochs in a 5000-epoch training cycle resulted in about a 25\% error reduction compared to the Adam-only optimizer at low and medium noise levels, demonstrating the efficacy of LBFGS during fine-tuning. However, this advantage diminished at higher noise levels, where the differences between the optimizers became less pronounced. In experiments where optimizers were switched periodically every 500 epochs, results varied, but overall, periodic switching helped minimize the error in certain cases, though it did not consistently outperform the Adam-LBFGS approach. Notably, using only LBFGS led to significantly higher errors, suggesting that while LBFGS is beneficial for fine-tuning, it tends to settle in local minima and does not further improve the error.

\subsection{Progressive roll-out length versus Constant roll-out length}
Figure~\cref{fig:abl_4} compares the performance of the model using progressive roll-out length scheduling and constant roll-out length scheduling under different noise levels. Progressive roll-out scheduling shows a 20\% improvement in mean training errors and a 16.67\% improvement in mean test errors compared to the constant roll-out length, indicating better learning and generalization. Although the time taken for training with progressive roll-out scheduling is 6.94\% lower, the overall training duration remains slightly higher due to the incremental increase in roll-out length. These results highlight the effectiveness of progressive roll-out scheduling in enhancing model performance.

\subsection{Overall Impact of Optimized Training Strategy}
\begin{figure}[!t]
\centering
\subfloat[Training Error]{\includegraphics[width=1\linewidth]{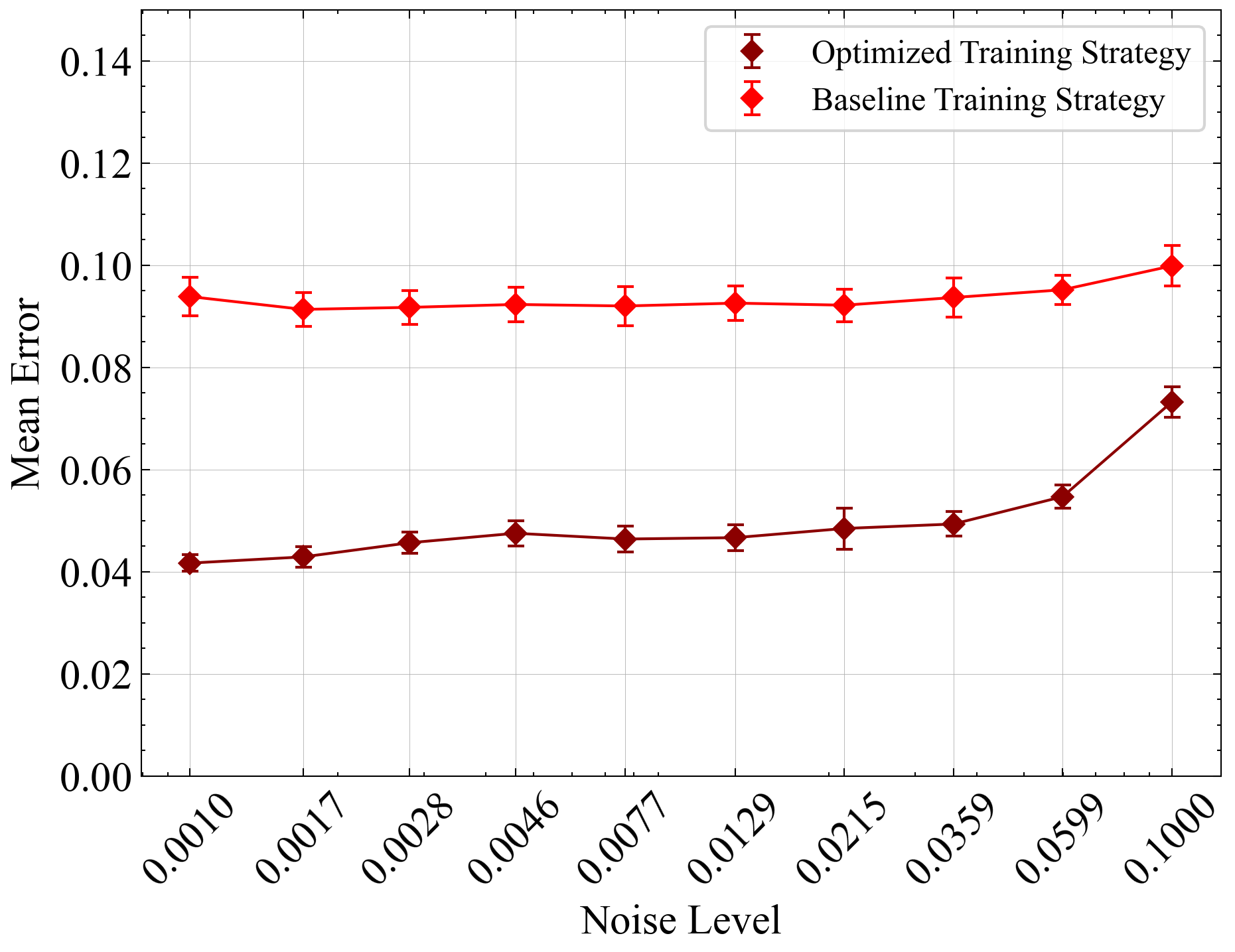}%
\label{fig:overall_train}}
\vfill
\subfloat[Test Error]{\includegraphics[width=1\linewidth]{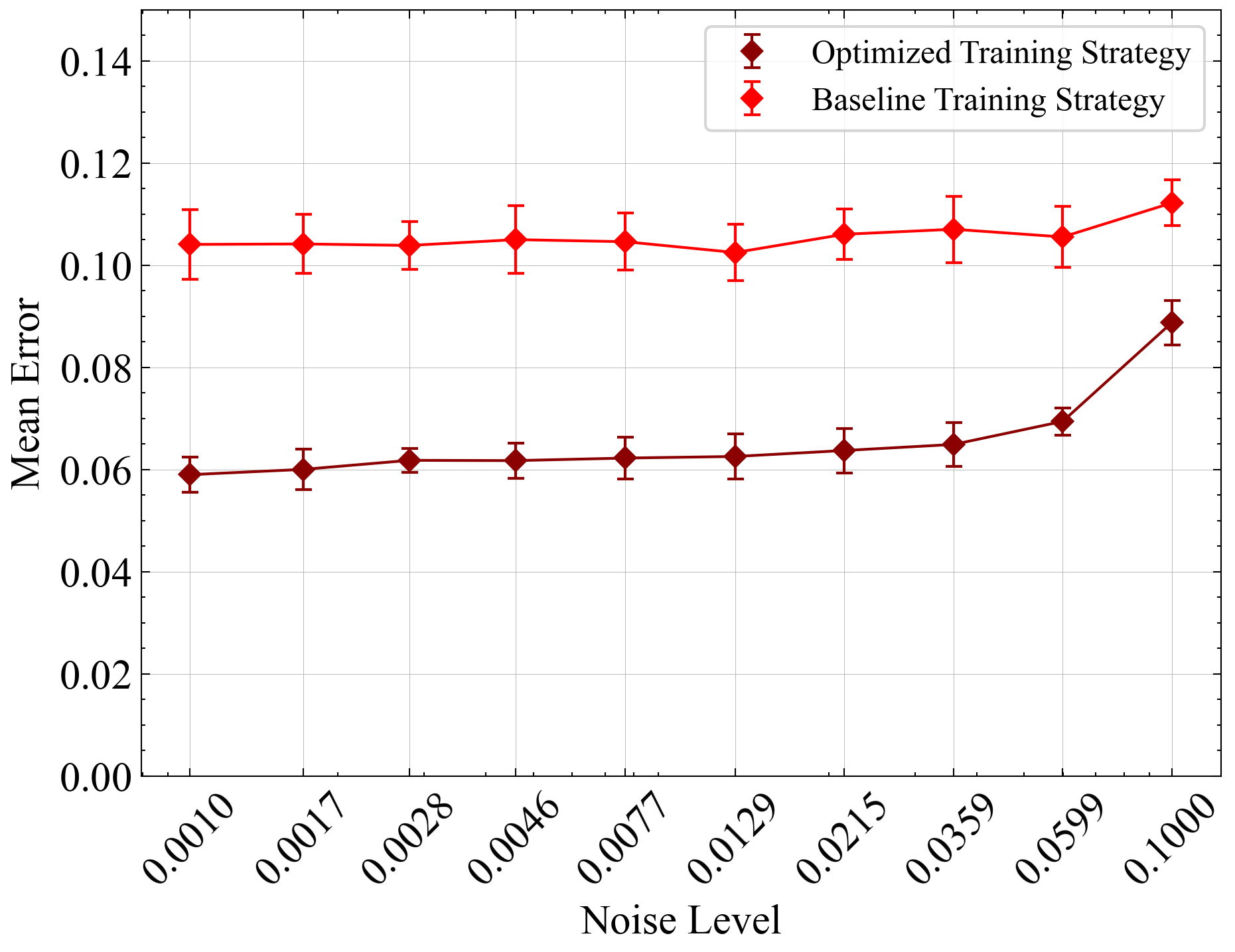}%
\label{fig:overall_test}}
\caption{Comparison of model prediction performance for the Van der Pol system. The “baseline” refers to the standard Koopman operator trained with fixed rollout length (R=99), no learning rate scheduling, and constant learning rate of 0.01. The “optimized” strategy incorporates a dissipative Koopman operator, progressive rollout, and cyclic learning rate schedule. The optimized configuration shows a 39\% improvement in training error and a 37.6\% improvement in testing error.}
\label{fig:overall}
\end{figure}

\begin{figure}[!t]
\centering
\includegraphics[width=0.9\linewidth]{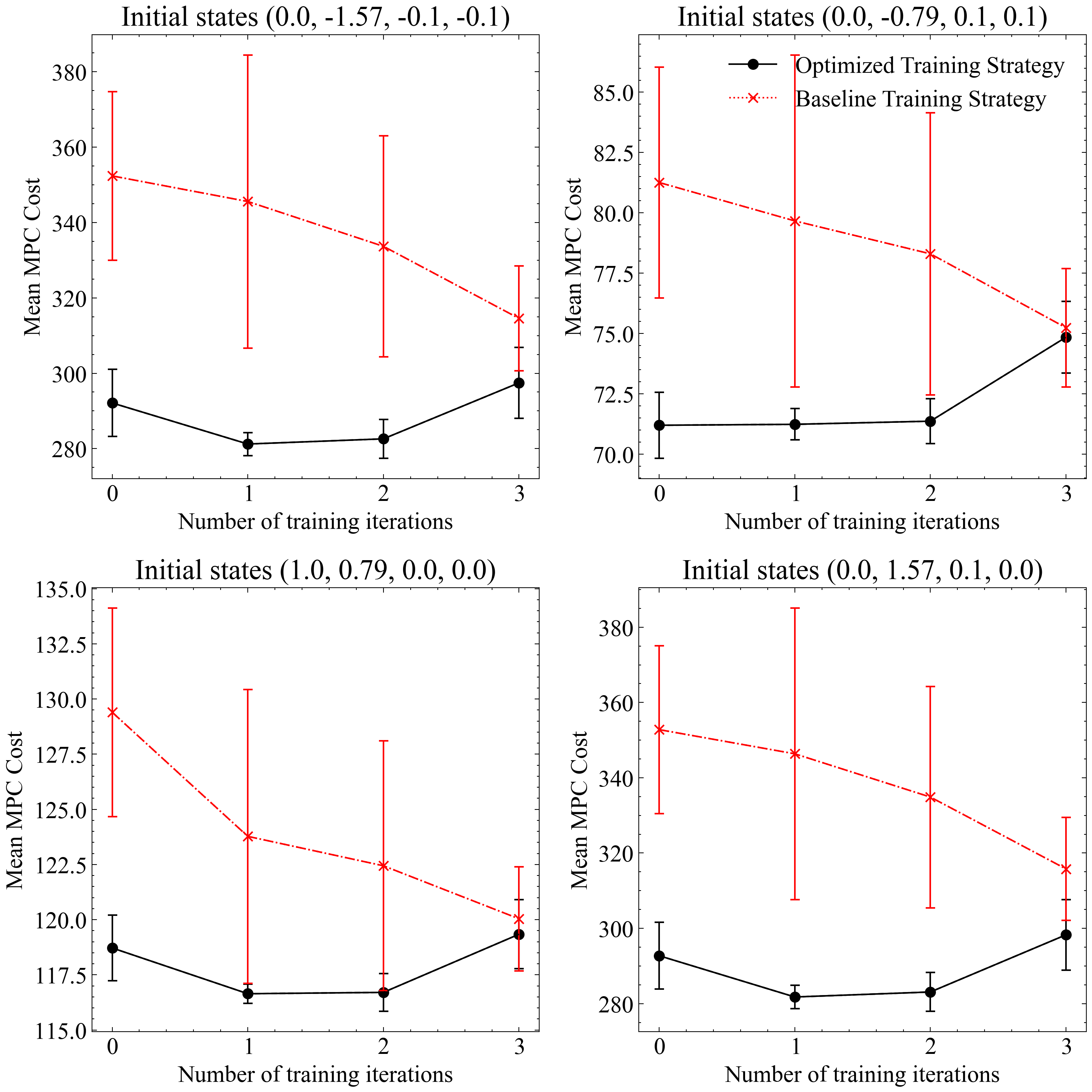}
\caption{Model predictive control performance for the Cartpole system. The “baseline” configuration uses the standard Koopman operator with fixed rollout length (R=16), constant learning rate, and no learning rate scheduling. The “optimized” configuration employs a dissipative Koopman operator, progressively increasing rollout length, a cyclic learning rate schedule, and fine-tuning with the LBFGS optimizer after initial training with Adam. Performance is evaluated across multiple iterative training data augmentations and initial conditions.}
\label{fig:overall_2}
\end{figure}

To assess the combined effect of the above-mentioned training strategies, we compare the overall performance of our optimized configuration with "baseline" configurations for Van der Pol system \cref{fig:overall} and Cartpole system \cref{fig:overall_2}. \Cref{fig:overall} shows the prediction performance comparison for the Van der Pol system and \cref{fig:overall_2} shows the model predictive costs comparison for the Cartpole system over iterative data augmentation. 

For the Van der Pol system \cref{fig:overall}, the baseline configuration utilized the \emph{Standard Koopman Operator} approach with a fixed maximum roll-out duration of 99, without any learning rate scheduler and maintained a constant learning rate of 0.01. In contrast, our optimized configuration employed a dissipative Koopman operator alongside a progressive roll-out strategy, starting with a roll-out length of 1 and geometrically increasing by a factor of 2 every 200 epochs until reaching the maximum roll-out length. Additionally, we implemented a triangular mode cyclic learning rate scheduler, with the maximum learning rate set to 0.1, halving every 500 epochs, and a base learning rate of 0.01. Both configurations used the Adam optimizer. In general, the optimized configuration resulted in a 39\% improvement in training error and a 37.6\% improvement in testing error compared to baseline. This significant reduction in error highlights how the combined components of the optimized configuration effectively enhance the model's performance. 

In the case of the Cartpole system (\cref{fig:overall_2}, the baseline configuration used the Adam optimizer with a constant roll-out duration of 16 and a constant learning rate. The \emph{Standard Koopman Operator} approach was used throughout the training. Conversely, the optimized strategy employed the Adam optimizer for the initial 1500 epochs, followed by fine-tuning using the LBFGS optimizer. The duration of the roll-out was progressively increased from 1 to 16, a cyclic learning rate was implemented, and a dissipative Koopman operator was used. The optimized training strategy outperforms the baseline training strategy in different initial conditions in multiple iterative training data augmentations.

\section{Conclusions}
\label{sec:Conclusions}
In this work, we proposed a novel data-driven learning framework for the Koopman operator of nonlinear dynamical systems. 
Unlike existing frameworks that rely on ad hoc observables or complex neural networks to approximate the Koopman invariant subspace, our approach utilizes observables generated by the system dynamics through a Hankel matrix. This data-driven time-delay embedding significantly enhances the performance and robustness of the resulting models compared to ad-hoc choices. When system dynamics is unknown, we approximate them using a neural network while benefiting from the structural feature of the Hankel matrix.
To ensure the robustness of the noise and the long-term stability of the lifted system, we designed a stable parameterization of the Koopman operator and employed a recurrent roll-out.
Our stabilized parameterization guarantees a stable system which is equally accurate or, in some cases, even more accurate than the unconstrained counterpart. 
Additionally, we addressed the curse of dimensionality in sampling from high-dimensional phase space by iteratively augmenting data from closed-loop trajectories. This approach not only improved performance in seen closed-loop scenarios, but also enhanced the model's ability to handle unseen control objectives effectively.
We conducted an ablation study to demonstrate the effectiveness of several key components: switching to the LBFGS optimizer for fine-tuning, employing cyclic learning rate scheduling, utilizing progressive roll-out lengths for enhanced training outcomes, and incorporating a stabilized Koopman operator to ensure guaranteed system stability.
Our numerical results on several systems with both known and unknown dynamics demonstrated the superior performance of the proposed framework in the prediction and model predictive control (MPC) of nonlinear systems, compared to other popular DMD-based methods. 
In summary, our framework offers a robust noise approach for learning stable Koopman operators, providing a solid foundation for future research (e.g., Koopman for reinforcement learning~\cite{rozwood2024koopman}) and applications in nonlinear dynamical system control. Future work may extend our approach to high-dimensional systems by integrating dimensionality reduction techniques (such as autoencoders or proper orthogonal decomposition) and leveraging parallelized computation, thereby broadening its applicability to real-world, large-scale problems.

\section*{Acknowledgement}
This work used NPL at RPI through allocation AGRM from the Center for Computational Innovations during the early stages of this research. Numerical experiments were performed on the FASTER and Anvil supercomputers at TAMU and Purdue, respectively, using computational resources granted by the Advanced Cyberinfrastructure Coordination Ecosystem: Services \& Support (ACCESS) program under allocation PHY240112 (supported by U.S. National Science Foundation grants \#2138259, \#2138286, \#2138307, \#2137603, and \#2138296) and by the National Energy Research Scientific Computing Center (NERSC), a DOE Office of Science User Facility, under award DDR-ERCAP0030714. The authors appreciate the use of the research computing resources of the Empire AI Consortium, Inc. with support from Empire State Development of the State of New York, the Simons Foundation, and the Secunda Family Foundation.

\bibliographystyle{IEEEtran}
\bibliography{panlab, new_citation}

\newpage

\section{Biography Section}

\begin{IEEEbiographynophoto}{Shahriar Akbar Sakib}
 is currently a second-year Ph.D. student in Mechanical
Engineering at Rensselaer Polytechnic Institute in Troy, NY. He
earned his B.Sc. in Mechanical Engineering from the Bangladesh University
of Engineering and Technology in 2021. His research focuses on data-driven
modeling and control of complex nonlinear dynamical systems, as well as
machine learning.
\end{IEEEbiographynophoto}

\begin{IEEEbiographynophoto}{Shaowu Pan}
received a B.S. degree in Aerospace Engineering and Applied
Mathematics from Beihang University, China, in 2013. He earned M.S. and
Ph.D. degrees in Aerospace Engineering and Scientific Computing from the
University of Michigan, Ann Arbor, in April 2021. Following this, he was
a Postdoctoral Scholar at the AI Institute in Dynamic Systems, University
of Washington, Seattle. His research interests include data-driven modeling
and control of complex nonlinear systems, scientific machine learning, and
dynamical systems.
\end{IEEEbiographynophoto}

\vfill

\end{document}